# On the Cauchy-Bunyakovsky-Schwarz Inequality in Pseudo-Euclidean Space

Ilia Lomidze[1]*, Natela Chachava[1]

[1]*School of Physics, Tbilisi Free University, Tbilisi, Georgia*

In a pseudo-Euclidean space, the Cauchy-Bunyakovsky-Schwarz inequality's generalization is derived and studied. Complex vector space is considered as well. Generalizations of angles and dihedral angles in these spaces are studied.

**Key words**: CBS inequality; pseudo-Euclidean space; degenerated space; complex space; angles in vector spaces; polyvectors; Lagrange identity; Gram's matrix.

**MSC(2020):** 15A03, 15A63, 15A75, 26D15, 30L15.

**1. Introduction.** It is well-known, that the Cauchy-Bunyakovsky-Schwarz (CBS) inequality

$$\left(\sum_{i=1}^{n} a_i b_i\right)^2 \le \left(\sum_{i=1}^{n} a_i^2\right)\left(\sum_{i=1}^{n} b_i^2\right) \tag{1}$$

is valid for any real numbers $a_1,\dots,a_n,b_1,\dots,b_n$ and the equality is achieved iff $a_i = kb_i,\ \forall k\in\mathbb{R},\ i=\overline{1,n}$. The geometric interpretation of inequality (1) considers these numbers $a_k, b_k,\ k=\overline{1,n},$ as Cartesian components of vectors $\boldsymbol{a}=\boldsymbol{e}_k a_k,\ \boldsymbol{b}=\boldsymbol{e}_k b_k$ (summation over repeating indices is implied), in Euclidean vector space $\mathbb{E}^n$. Moreover, (1) provides the opportunity to define the angle between vectors $\boldsymbol{a},\boldsymbol{b}\in\mathbb{E}^n$:

$$\cos\varphi = a_i b_i \left(a_j a_j\right)^{-1/2}\left(b_k b_k\right)^{-1/2}. \tag{2}$$

We have studied relations analogous to the inequality (1) for vectors $\boldsymbol{a},\boldsymbol{b}\in\mathbb{P}^{n+m}_{n,m}$, where $\mathbb{P}^{n+m}_{n,m},\ n,m\in\mathbb{N}$, denotes a pseudo-Euclidean vector space with a Cartesian metric

$$g_{ij} = g^{ij} = \mathrm{diag}(\overbrace{+1,\dots,+1}^{n},\overbrace{-1,\dots,-1}^{m}),\quad i,j=\overline{1,n+m}. \tag{3}$$

Bellow we use notations $\langle \boldsymbol{a}\,|\,\boldsymbol{b}\rangle$ for an inner product of any (complex) vectors $\boldsymbol{a},\boldsymbol{b}\in\mathbb{P}^{n+m}_{n,m}$:

$$\langle \boldsymbol{a}|\boldsymbol{b}\rangle = \bar{a}^j b_j = \overline{\langle \boldsymbol{b}|\boldsymbol{a}\rangle},\quad \boldsymbol{a}^2 = \langle \boldsymbol{a}\,|\,\boldsymbol{a}\rangle\in\mathbb{R}, \tag{4}$$

$$a_i = \langle \boldsymbol{e}_i|\boldsymbol{a}\rangle = \langle \boldsymbol{e}_i|\boldsymbol{e}_j\rangle a^j = g_{ij}a^j,\quad a^j = g^{ij}a_i,\quad (a_i, b_i\in\mathbb{C}) \tag{4'}$$

where bar over a symbol means its complex conjugate and $\{\boldsymbol{e}_i\,|\,\boldsymbol{e}_i\in\mathbb{P}^{n+m}_{n,m},\ i=\overline{1,n+m}\}$ is a basis in $\mathbb{P}^{n+m}_{n,m}$. In a real vector space, one obviously has

$$\bar{a}_i = a_i,\ \langle \boldsymbol{a}|\boldsymbol{b}\rangle = \langle \boldsymbol{b}|\boldsymbol{a}\rangle.$$

Vectors having positive (negative) squares we call *time-like* and denote as $\boldsymbol{t}$ (*space-like* and denote as $\boldsymbol{s}$), as it is assumed in Physics; vectors $\boldsymbol{c}\neq\boldsymbol{0}$, such that $\boldsymbol{c}^2 = \langle \boldsymbol{c}\,|\,\boldsymbol{c}\rangle^2 = 0$ we call *isotropic* or *light-like.*

For a norm of the vector $\boldsymbol{a}\in\mathbb{P}^{n+m}_{n,m}$ we use a notation $|\boldsymbol{a}|$. So, in pseudo-Euclidean space $\mathbb{P}^{n+m}_{n,m}$

$$|\boldsymbol{a}| = \sqrt{\boldsymbol{a}^2},\ \text{if } \boldsymbol{a} \text{ is time-like } (\boldsymbol{a}^2>0);\quad |\boldsymbol{a}| = \sqrt{-\boldsymbol{a}^2},\ \text{if } \boldsymbol{a} \text{ is space-like } (\boldsymbol{a}^2<0);\quad |\boldsymbol{c}| = \sqrt{\boldsymbol{c}^2} = 0.$$

**Note.** (Complex) 2-plane is pseudo-Euclidean iff it contains two linearly independent isotropic vectors

$$\boldsymbol{c}_j = \boldsymbol{t} + \boldsymbol{s}e^{i\varphi_j},\ \forall\varphi_j\in(-\pi,\pi],\ \langle \boldsymbol{t}|\boldsymbol{s}\rangle = 0,\ \boldsymbol{t}^2 = -\boldsymbol{s}^2 = 1,\ \boldsymbol{c}_j^2 = 0,\ j=1,2,\ \boldsymbol{c}_1\neq k\boldsymbol{c}_2,\ k\in\mathbb{C},$$

and no linear combinations of them are mutually orthogonal to any other vector of the plane:

$$\langle \boldsymbol{c}\,|\,\alpha\boldsymbol{t}+\beta\boldsymbol{s}\rangle = \alpha\langle \boldsymbol{t}|\boldsymbol{t}\rangle + \beta e^{-i\varphi}\langle \boldsymbol{s}|\boldsymbol{s}\rangle = \alpha - \beta e^{-i\varphi} = 0 \Leftrightarrow \beta = \alpha e^{i\varphi} \Leftrightarrow \alpha\boldsymbol{t}+\beta\boldsymbol{s} = \alpha\boldsymbol{c},\ \forall\alpha\in\mathbb{C},$$

$$\langle \boldsymbol{c}_1|\boldsymbol{c}_2\rangle = \langle \boldsymbol{t}|\boldsymbol{t}\rangle + \langle \boldsymbol{s}|\boldsymbol{s}\rangle e^{i(\varphi_2-\varphi_1)} = 1 - e^{i(\varphi_2-\varphi_1)} \neq 0\ (\varphi_2\neq\varphi_1).$$

*E-mail: i.lomidze@freeuni.edu.ge

It is well-known number of generalizations of CSB inequality (see, e. g., [15]) for real Euclidean space as well as for a complex one. Some authors consider inequalities for sets of vectors those contain more then two elements, others – for sequences those do not being vectors of any linear space. In [1]-[4] the Minkowski space (pseudo-Euclidean space $\mathbb{P}^4_{1,3}$, widely used in Physics) is considered. It is asserted (and proven) therein that for time-like vectors, the CBS inequality reverses its sign.

However, to the best of our knowledge, no comprehensive study of CBS-type relations for vectors in an arbitrary pseudo-Euclidean space has been conducted to date. In this paper, we address this gap by investigating such relations in real and complex pseudo-Euclidean spaces, as well as in their degenerate subspaces. Based on these results, we introduce and examine the concept of angles (in our view, the definitions of the Lorentzian time-like angle introduced in [4] and [5] require clarification) between arbitrary vectors in a pseudo-Euclidean space. Our proposed approach also opens up new possibilities for the study of multi-vector spaces – real and complex Euclidean, as well as pseudo-Euclidean ones.

**2. Main results.** We established that in real pseudo-Euclidean space $\mathbb{P}^{n+m}_{n,m}$ is valid the following

**Theorem 1.** 1) *For vectors* $\boldsymbol{a},\boldsymbol{b}$ ($\boldsymbol{a}\neq k\boldsymbol{b}$, $k\in\mathbb{R}$) *the inequality* $\langle \boldsymbol{a}\,|\,\boldsymbol{b}\rangle^2<\boldsymbol{a}^2\boldsymbol{b}^2$ *holds iff* $\boldsymbol{a},\boldsymbol{b}\in\mathbb{P}^2_{0,2}\subset\mathbb{P}^{n+m}_{n,m}$ *or* $\boldsymbol{a},\boldsymbol{b}\in\mathbb{P}^2_{2,0}\subset\mathbb{P}^{n+m}_{n,m}$; 2) *the inequality* $\langle \boldsymbol{a}\,|\,\boldsymbol{b}\rangle^2>\boldsymbol{a}^2\boldsymbol{b}^2$ *holds iff* $\boldsymbol{a},\boldsymbol{b}\in\mathbb{P}^2_{1,1}\subset\mathbb{P}^{n+m}_{n,m}$ ($\boldsymbol{a}\neq k\boldsymbol{b}$); 3) *the equality* $\langle \boldsymbol{a}\,|\,\boldsymbol{b}\rangle^2=\boldsymbol{a}^2\boldsymbol{b}^2$ *holds iff* $\boldsymbol{a},\boldsymbol{b}$ ($\boldsymbol{a}\neq k\boldsymbol{b}$) *belong to* $\operatorname{span}[\boldsymbol{s},\boldsymbol{c}]$, $\langle \boldsymbol{c}\,|\,\boldsymbol{s}\rangle=0$ *or to* $\operatorname{span}[\boldsymbol{t},\boldsymbol{c}]$, $\langle \boldsymbol{c}\,|\,\boldsymbol{t}\rangle=0$ *or to* $\operatorname{span}[\boldsymbol{c}_1,\boldsymbol{c}_2]$, $\langle \boldsymbol{c}_1\,|\,\boldsymbol{c}_2\rangle=0$, $\boldsymbol{c}_1^2=\boldsymbol{c}_2^2=0$ (***s-c*** *kind,* ***t-c*** *kind or* ***c-c*** *kind* 2*-planes, respectively*). 4) *For any vectors* $\boldsymbol{a}=k\boldsymbol{b}\in\mathbb{P}^{n+m}_{n,m}$ ($k\in\mathbb{R}$) *the equality* $\langle \boldsymbol{a}\,|\,\boldsymbol{b}\rangle^2=\boldsymbol{a}^2\boldsymbol{b}^2$ *holds.*

**Note.** All such 2-planes do exist in pseudo-Euclidean space $\mathbb{P}^{n+m}_{n,m}$ (see, e.g., [1]).

**Proof** (real case). Considering $\{\boldsymbol{u}\,|\,\boldsymbol{u}=\alpha\boldsymbol{a}+\beta\boldsymbol{b},\ \boldsymbol{a}\neq k\boldsymbol{b},\ k,\alpha,\beta\in\mathbb{R},\ \alpha^2+\beta^2\neq 0\}\subset\mathbb{P}^{n+m}_{n,m}$, in order to determine the sign of $\boldsymbol{u}^2=\alpha^2\boldsymbol{a}^2+2\alpha\beta\langle \boldsymbol{a}|\boldsymbol{b}\rangle+\beta^2\boldsymbol{b}^2$, one has to consider the following special cases:

1) $D\equiv\langle \boldsymbol{a}|\boldsymbol{b}\rangle^2-\boldsymbol{a}^2\boldsymbol{b}^2<0 \ \Leftrightarrow\ \boldsymbol{a}^2\boldsymbol{b}^2>\langle \boldsymbol{a}|\boldsymbol{b}\rangle^2\geq 0$. Then $\boldsymbol{u}^2>0$ if $\boldsymbol{a}^2,\boldsymbol{b}^2>0$ and $\boldsymbol{u}^2<0$ if $\boldsymbol{a}^2,\boldsymbol{b}^2<0$. The inverse is obvious, namely, if $\boldsymbol{a},\boldsymbol{b}\in\mathbb{P}^2_{2,0}$ (or $\boldsymbol{a},\boldsymbol{b}\in\mathbb{P}^2_{0,2}$), then $\langle \boldsymbol{a}|\boldsymbol{b}\rangle^2-\boldsymbol{a}^2\boldsymbol{b}^2\leq 0$ and the equality achieves iff $\boldsymbol{a}=k\boldsymbol{b}$, $k\in\mathbb{R}$, just as the Theorem states.

2) $\langle \boldsymbol{a}|\boldsymbol{b}\rangle^2>\boldsymbol{a}^2\boldsymbol{b}^2$. Then $D=\langle \boldsymbol{a}|\boldsymbol{b}\rangle^2-\boldsymbol{a}^2\boldsymbol{b}^2>0$ and the $\operatorname{span}[\boldsymbol{a},\boldsymbol{b}]$ contains vectors of all three types. Therefore, in this case the $\operatorname{span}[\boldsymbol{a},\boldsymbol{b}]$ is pseudo-Euclidean 2-plane. Inversely, if the $\operatorname{span}[\boldsymbol{a},\boldsymbol{b}]$ is pseudo-Euclidean 2-plane, that is, if

$$\boldsymbol{a}=\lambda\boldsymbol{t}+\rho\boldsymbol{s},\ \boldsymbol{b}=\mu\boldsymbol{t}+\sigma\boldsymbol{s},\ \boldsymbol{t}^2=1,\ \boldsymbol{s}^2=-1,\ \langle \boldsymbol{t}|\boldsymbol{s}\rangle=0,\ \ (\lambda,\mu,\rho,\sigma\in\mathbb{R})$$

then

$$\langle \boldsymbol{a}|\boldsymbol{b}\rangle^2-\boldsymbol{a}^2\boldsymbol{b}^2=-2\lambda\mu\rho\sigma+\lambda^2\sigma^2+\rho^2\mu^2=(\lambda\sigma-\rho\mu)^2\geq 0,$$

and the equality achieves iff

$$\lambda\sigma=\rho\mu\ \Leftrightarrow\ \boldsymbol{a}=k\boldsymbol{b},\ \ k=\lambda/\mu=\rho/\sigma\in\mathbb{R},$$

just as the Theorem states.

**Note.** In the case under consideration the quadratic equation for $\boldsymbol{u}\in\operatorname{span}[\boldsymbol{a},\boldsymbol{b}]$, $\boldsymbol{u}^2=0$, that is

$$\alpha^2\boldsymbol{a}^2+2\alpha\beta\langle \boldsymbol{a}|\boldsymbol{b}\rangle+\beta^2\boldsymbol{b}^2=0 \tag{5}$$

has exactly 2 independent real solutions for coefficients $\alpha,\beta\in\mathbb{R}$. Correspondingly, 2-plane $\operatorname{span}[\boldsymbol{a},\boldsymbol{b}]$ contains exactly 2 linearly independent isotropic vectors that are not mutually orthogonal, just as it must be for real pseudo-Euclidean 2-plane.

3) $\langle \boldsymbol{a}|\boldsymbol{b}\rangle^2 = \boldsymbol{a}^2\boldsymbol{b}^2 \ge 0 \;\Leftrightarrow\; D = \langle \boldsymbol{a}|\boldsymbol{b}\rangle^2 - \boldsymbol{a}^2\boldsymbol{b}^2 = 0,\; \boldsymbol{a} \ne k\boldsymbol{b},\; k \in \mathbb{R}.$

$3^a$) If $\boldsymbol{b}^2 = 0$, then

$$\langle \boldsymbol{a}|\boldsymbol{b}\rangle = 0 \;\Rightarrow\; \langle \alpha\boldsymbol{a}+\beta\boldsymbol{b}|\boldsymbol{b}\rangle = 0,\;\; (\alpha\boldsymbol{a}+\beta\boldsymbol{b})^2 = \alpha^2\boldsymbol{a}^2,\;\; \forall \alpha,\beta \in \mathbb{R}.$$

Hence, in this case the type of the vector $\boldsymbol{u} \in \mathrm{span}[\boldsymbol{a},\boldsymbol{b}]$ (the sign of $\boldsymbol{u}^2$) is completely determined by the sign of $\boldsymbol{a}^2$: the $\mathrm{span}[\boldsymbol{a},\boldsymbol{b}]$ (2-plane) consists of all vectors $\boldsymbol{u}$ those are orthogonal to the isotropic vector $\boldsymbol{b}$ (unique in this plane), all $\boldsymbol{u}$ being time-like (***c-t*** kind degenerate 2-plane) if time-like is the vector $\boldsymbol{a}$ or all $\boldsymbol{u}$ being space-like (***c-s*** kind degenerate 2-plane) if space-like is the vector $\boldsymbol{a}$.

If the vector $\boldsymbol{a}$ itself is isotropic, then all vectors $\boldsymbol{u} \in \mathrm{span}[\boldsymbol{a},\boldsymbol{b}]$ are mutually orthogonal isotropic ones, hence, the correspondent 2-plane is isotropic (***c-c*** kind degenerate 2-plane).

$3^b$) If $\boldsymbol{b}^2 \ne 0$, then, due to the condition $D = \langle \boldsymbol{a}|\boldsymbol{b}\rangle^2 - \boldsymbol{a}^2\boldsymbol{b}^2 = 0$, the equation (5) has double (unique) root $\beta_0 = -\alpha\langle \boldsymbol{a}|\boldsymbol{b}\rangle/\boldsymbol{b}^2$ and hence

$$\boldsymbol{u}(\beta_0) = \alpha\left(\boldsymbol{a} - \boldsymbol{b}\langle \boldsymbol{a}|\boldsymbol{b}\rangle/\boldsymbol{b}^2\right) \;\Rightarrow\; \left\{\langle \boldsymbol{u}(\beta_0)|\boldsymbol{a}\rangle = 0,\; \langle \boldsymbol{u}(\beta_0)|\boldsymbol{b}\rangle = 0\right\},$$

that is, any vector of the $\mathrm{span}[\boldsymbol{a},\boldsymbol{b}]$ is orthogonal to the unique isotropic vector $\boldsymbol{u}(\beta_0)$. Besides, quadratic trinomial – the left-hand side (LHS) of the eq. (5) retains positive (negative) sign iff the vectors $\boldsymbol{a}$ and $\boldsymbol{b}$ both are ***t***-like (***s***-like). So, in this case the $\mathrm{span}[\boldsymbol{a},\boldsymbol{b}]$ is a ***c-t*** kind (***c-s*** kind) degenerate 2-plane.

Inversely, if the $\mathrm{span}[\boldsymbol{a},\boldsymbol{b}]$ is a ***c-t*** kind (***c-s*** kind) degenerate 2-plane, that is, if

$$\boldsymbol{a} = \lambda\boldsymbol{t} + \rho\boldsymbol{c},\quad \boldsymbol{b} = \mu\boldsymbol{t} + \sigma\boldsymbol{c} \quad (\text{or } \boldsymbol{a} = \lambda\boldsymbol{s} + \rho\boldsymbol{c},\; \boldsymbol{b} = \mu\boldsymbol{s} + \sigma\boldsymbol{c}),$$

$$\boldsymbol{t}^2 = 1,\; \boldsymbol{s}^2 = -1,\; \boldsymbol{c}^2 = 0,\; \langle \boldsymbol{t}|\boldsymbol{c}\rangle = 0,\; \langle \boldsymbol{s}|\boldsymbol{c}\rangle = 0\;\; (\lambda,\mu,\rho,\sigma \in \mathbb{R}),$$

then

$$\langle \boldsymbol{a}|\boldsymbol{b}\rangle^2 - \boldsymbol{a}^2\boldsymbol{b}^2 = (\pm\lambda\mu)^2 - (\pm\lambda^2)(\pm\mu^2) = 0. \qquad \square$$

**Theorem 2.** In complex pseudo-Euclidean vector space $\mathbb{P}^{n+m}_{n,m}$ the following is valid:

1) *For vectors* $\boldsymbol{a},\boldsymbol{b}$ ($\boldsymbol{a} \ne k\boldsymbol{b},\, k \in \mathbb{C}$) *the inequality* $|\langle \boldsymbol{a}|\boldsymbol{b}\rangle|^2 < \boldsymbol{a}^2\boldsymbol{b}^2$ *holds iff* $\boldsymbol{a},\boldsymbol{b} \in \mathbb{P}^2_{0,2} \subset \mathbb{P}^{n+m}_{n,m}$ *or* $\boldsymbol{a},\boldsymbol{b} \in \mathbb{P}^2_{2,0} \subset \mathbb{P}^{n+m}_{n,m}$; 2) *the inequality* $|\langle \boldsymbol{a}|\boldsymbol{b}\rangle|^2 > \boldsymbol{a}^2\boldsymbol{b}^2$ *holds iff* $\boldsymbol{a},\boldsymbol{b} \in \mathbb{P}^2_{1,1} \subset \mathbb{P}^{n+m}_{n,m}$ ($\boldsymbol{a} \ne k\boldsymbol{b}$); 3) *the equality* $|\langle \boldsymbol{a}|\boldsymbol{b}\rangle|^2 = \boldsymbol{a}^2\boldsymbol{b}^2$ *holds iff* $\boldsymbol{a},\boldsymbol{b}$ ($\boldsymbol{a} \ne k\boldsymbol{b}$) *belong to* $\mathrm{span}[\boldsymbol{s},\boldsymbol{c}]$, $\langle \boldsymbol{c}|\boldsymbol{s}\rangle = 0$ *or to* $\mathrm{span}[\boldsymbol{t},\boldsymbol{c}]$, $\langle \boldsymbol{c}|\boldsymbol{t}\rangle = 0$ *or to* $\mathrm{span}[\boldsymbol{c}_1,\boldsymbol{c}_2]$, $\langle \boldsymbol{c}_1|\boldsymbol{c}_2\rangle = 0$, $\boldsymbol{c}_1^2 = \boldsymbol{c}_2^2 = 0$, $\boldsymbol{c}_1 \ne k\boldsymbol{c}_2$ *(**s-c** kind, **t-c** kind or **c-c** kind 2-planes, respectively). In the last case in addition one gets* $|\langle \boldsymbol{a}|\boldsymbol{b}\rangle|^2 = \boldsymbol{a}^2\boldsymbol{b}^2 = 0$. 4) *For any vectors* $\boldsymbol{a} = k\boldsymbol{b} \in \mathbb{P}^{n+m}_{n,m}$ ($k \in \mathbb{C}$) *the equality* $|\langle \boldsymbol{a}|\boldsymbol{b}\rangle|^2 = \boldsymbol{a}^2\boldsymbol{b}^2$ *holds.*

Note that the Theorem 2 is very similar to the Theorem 1. The only difference is due to the definition of inner product of complex vectors that leads to the replacement $\langle \boldsymbol{a}|\boldsymbol{b}\rangle^2 \to |\langle \boldsymbol{a}|\boldsymbol{b}\rangle|^2$.

**Proof.** Considering again $\left\{\boldsymbol{u}\,|\,\boldsymbol{u} = \alpha\boldsymbol{a} + \beta\boldsymbol{b},\; \boldsymbol{a} \ne k\boldsymbol{b},\; k,\alpha,\beta \in \mathbb{C},\; |\alpha|^2 + |\beta|^2 \ne 0\right\} \subset \mathbb{P}^{n+m}_{n,m}$ and calculating $\boldsymbol{u}^2 = \langle \boldsymbol{u}|\boldsymbol{u}\rangle$ according to the definition of complex vectors' inner product, obviously, one gets (compare with the LHS of the equation (5)):

$$\begin{aligned}\boldsymbol{u}^2 = \langle \boldsymbol{u}|\boldsymbol{u}\rangle &= |\alpha|^2\langle \boldsymbol{a}|\boldsymbol{a}\rangle + |\beta|^2\langle \boldsymbol{b}|\boldsymbol{b}\rangle + \bar{\alpha}\beta\langle \boldsymbol{a}|\boldsymbol{b}\rangle + \alpha\bar{\beta}\langle \boldsymbol{b}|\boldsymbol{a}\rangle = \\ &= |\alpha|^2\left\{\boldsymbol{a}^2 + |\beta/\alpha|^2\boldsymbol{b}^2 + (\beta/\alpha)\langle \boldsymbol{a}|\boldsymbol{b}\rangle + \overline{(\beta/\alpha)\langle \boldsymbol{a}|\boldsymbol{b}\rangle}\right\}.\end{aligned} \tag{6}$$

Let firstly consider the case $\langle \boldsymbol{b}|\boldsymbol{b}\rangle = \boldsymbol{b}^2 \neq 0$. Then, applying the notation $z = \beta/\alpha$, we obtain

$$\begin{aligned} \boldsymbol{u}^2(z) &= |\alpha|^2 \boldsymbol{b}^2 \left[ \boldsymbol{a}^2/\boldsymbol{b}^2 + |z|^2 + 2\mathrm{Re}\left(z\langle \boldsymbol{a}|\boldsymbol{b}\rangle/\boldsymbol{b}^2\right) \right] = \\ &= |\alpha|^2 \boldsymbol{b}^2 \left[ |z + \langle \boldsymbol{b}|\boldsymbol{a}\rangle/\boldsymbol{b}^2|^2 - |\langle \boldsymbol{b}|\boldsymbol{a}\rangle/\boldsymbol{b}^2|^2 + \boldsymbol{a}^2/\boldsymbol{b}^2 \right] = \\ &= |\alpha|^2 \boldsymbol{b}^2 \left[ |z - z_0|^2 - R^2 \right]. \end{aligned} \tag{7}$$

In the last expression we denote

$$z_0 = -\langle \boldsymbol{b}|\boldsymbol{a}\rangle/\boldsymbol{b}^2, \tag{8}$$

$$R^2 = |\langle \boldsymbol{b}|\boldsymbol{a}\rangle/\boldsymbol{b}^2|^2 - \boldsymbol{a}^2/\boldsymbol{b}^2 = \left( |\langle \boldsymbol{a}|\boldsymbol{b}\rangle|^2 - \boldsymbol{a}^2\boldsymbol{b}^2 \right)/\langle \boldsymbol{b}|\boldsymbol{b}\rangle^2. \tag{9}$$

If the real

$$|\langle \boldsymbol{a}|\boldsymbol{b}\rangle|^2 - \boldsymbol{a}^2\boldsymbol{b}^2 \equiv D \tag{10}$$

is negative (hence, $R^2$ in (9) is negative, $R$ itself being purely imagine), then the expression in square brackets in (7) is positive for all $z \in \mathbb{C}$ and the sign of $\boldsymbol{u}^2(z)$ coincides with the sign of $\boldsymbol{b}^2 \neq 0$. Therefore, if $|\langle \boldsymbol{a}|\boldsymbol{b}\rangle|^2 < \boldsymbol{a}^2\boldsymbol{b}^2$, then the span$[\boldsymbol{a},\boldsymbol{b}]$ is Euclidean (complex) 2-plane (***t-t*** kind or ***s-s*** kind), similarly to the case of real vector space.

If $D > 0$, $\boldsymbol{b}^2 \neq 0$, then the sign of square brackets in the expression (7) is negative (positive) for $z$ inside (outside) the circle

$$|z - z_0| = R = \sqrt{D}/|\boldsymbol{b}^2|, \quad z \in \mathbb{C}. \tag{11}$$

For $z = z_0 + R\exp(i\varphi)$, $\forall \varphi \in \mathbb{R}$, (continuous set of such $z$) one gets $\boldsymbol{u}^2(z) = 0$, so we have continuous set of isotropic vectors $\boldsymbol{u}(z) \neq \boldsymbol{0}$ if $D > 0$ (it is easy to check that only 2 vectors of this set are linearly independent but not mutually orthogonal, see Appendix 1). That is, in this case the span$[\boldsymbol{a},\boldsymbol{b}]$ contains vectors of all kinds (***t***-like, ***s***-like and ***c***-like) and, i.e., it is pseudo-Euclidean (complex) 2-plane.

If $D = 0 \;\Leftrightarrow\; |\langle \boldsymbol{a}|\boldsymbol{b}\rangle|^2 = \boldsymbol{a}^2\boldsymbol{b}^2 \; (\boldsymbol{b}^2 \neq 0)$, then

$$\boldsymbol{u}^2(z)|_{z=z_0} = 0, \quad z_0 = -\langle \boldsymbol{b}|\boldsymbol{a}\rangle/\boldsymbol{b}^2.$$

In this case we have

$$\boldsymbol{u}(z_0) = \alpha(\boldsymbol{a} + z_0\boldsymbol{b}) = \alpha(\boldsymbol{a} - \boldsymbol{b}\langle \boldsymbol{b}|\boldsymbol{a}\rangle/\boldsymbol{b}^2) \qquad \Rightarrow$$

$$\langle \boldsymbol{u}(z_0)|\boldsymbol{a}\rangle = \bar{\alpha}\left[ \langle \boldsymbol{a}|\boldsymbol{a}\rangle - \left( \langle \boldsymbol{a}|\boldsymbol{b}\rangle/\boldsymbol{b}^2 \right)\langle \boldsymbol{b}|\boldsymbol{a}\rangle \right] = \bar{\alpha}\left[ \boldsymbol{a}^2 - |\langle \boldsymbol{a}|\boldsymbol{b}\rangle|^2/\boldsymbol{b}^2 \right] = -\bar{\alpha}D/\boldsymbol{b}^2 = 0,$$

$$\langle \boldsymbol{u}(z_0)|\boldsymbol{b}\rangle = \bar{\alpha}\left[ \langle \boldsymbol{a}|\boldsymbol{b}\rangle - \left( \langle \boldsymbol{a}|\boldsymbol{b}\rangle/\boldsymbol{b}^2 \right)\boldsymbol{b}^2 \right] = \bar{\alpha}\left[ \langle \boldsymbol{a}|\boldsymbol{b}\rangle - \langle \boldsymbol{a}|\boldsymbol{b}\rangle \right] = 0.$$

Hence, in the case $|\langle \boldsymbol{a}|\boldsymbol{b}\rangle|^2 = \boldsymbol{a}^2\boldsymbol{b}^2$ the (unique) isotropic vector $\boldsymbol{u}(z_0)$ is orthogonal to any vector of the span$[\boldsymbol{a},\boldsymbol{b}]$, that is, in this case the span$[\boldsymbol{a},\boldsymbol{b}]$ is ***c-b*** kind degenerate (complex) 2-plane.

Now let $\langle \boldsymbol{b}|\boldsymbol{b}\rangle = \boldsymbol{b}^2 = 0$, $\boldsymbol{b} \neq \boldsymbol{0}$. Then, according to (10), $D = |\langle \boldsymbol{a}|\boldsymbol{b}\rangle|^2 \geq 0$ and

$$\boldsymbol{u}^2 = |\alpha|^2 \left\{ \boldsymbol{a}^2 + 2\mathrm{Re}\left( z\langle \boldsymbol{a}|\boldsymbol{b}\rangle \right) \right\}. \tag{12}$$

If $D = 0 \;\Leftrightarrow\; \langle \boldsymbol{a}|\boldsymbol{b}\rangle = 0$ and therefore

$$\langle \boldsymbol{u}(z)|\boldsymbol{a}\rangle = \bar{\alpha}\left[ \boldsymbol{a}^2 + z\langle \boldsymbol{b}|\boldsymbol{a}\rangle \right] = \bar{\alpha}\boldsymbol{a}^2, \; \langle \boldsymbol{u}(z)|\boldsymbol{b}\rangle = \bar{\alpha}\left[ \langle \boldsymbol{a}|\boldsymbol{b}\rangle + z\langle \boldsymbol{b}|\boldsymbol{b}\rangle \right] = 0, \; \boldsymbol{u}^2(z) = |\alpha|^2\boldsymbol{a}^2, \; \forall z \in \mathbb{C},$$

hence, for $D = 0$, $\boldsymbol{b}^2 = 0$, $\boldsymbol{b} \neq \boldsymbol{0}$, the span$[\boldsymbol{a},\boldsymbol{b}]$ is ***c-a*** kind degenerate (complex) 2-plane.

If $D>0,\ \boldsymbol{b}^2=0\ (\boldsymbol{b}\neq\boldsymbol{0}) \Rightarrow \langle\boldsymbol{a}|\boldsymbol{b}\rangle\neq 0$ and hence for $\forall\operatorname{Im}(z\langle\boldsymbol{a}|\boldsymbol{b}\rangle)$ we get from (12)

$$\operatorname{Re}(z\langle\boldsymbol{a}|\boldsymbol{b}\rangle) > -\boldsymbol{a}^2/2 \ \Leftrightarrow\ \boldsymbol{u}^2>0,$$
$$\operatorname{Re}(z\langle\boldsymbol{a}|\boldsymbol{b}\rangle) < -\boldsymbol{a}^2/2 \ \Leftrightarrow\ \boldsymbol{u}^2<0,$$
$$\operatorname{Re}(z\langle\boldsymbol{a}|\boldsymbol{b}\rangle) = -\boldsymbol{a}^2/2 \ \Leftrightarrow\ \boldsymbol{u}^2=0.$$

Continuous set of isotropic vectors

$$\boldsymbol{u}(z)=\alpha(\boldsymbol{a}+z\boldsymbol{b}),\ \ \boldsymbol{b}^2=0,\ \ \operatorname{Re}(z\langle\boldsymbol{a}|\boldsymbol{b}\rangle)=-\boldsymbol{a}^2/2,\quad \forall\operatorname{Im}(z\langle\boldsymbol{a}|\boldsymbol{b}\rangle),$$

obviously, contains 2 linearly independent (complex) vectors that are not mutually orthogonal. That is, in the case under consideration the span$[\boldsymbol{a},\boldsymbol{b}]$ is a pseudo-Euclidean (complex) 2-plane.

So, we have proved that for any complex vectors $\boldsymbol{a},\boldsymbol{b}\in\mathbb{P}^{n+m}_{n,m}$, $\boldsymbol{a}\neq k\boldsymbol{b}$, $k\in\mathbb{C}$, the span$[\boldsymbol{a},\boldsymbol{b}]$ is:

a pseudo-Euclidean (complex) 2-plane, if

$$D=|\langle\boldsymbol{a}|\boldsymbol{b}\rangle|^2-\boldsymbol{a}^2\boldsymbol{b}^2>0 \ \Leftrightarrow\ |\langle\boldsymbol{a}|\boldsymbol{b}\rangle|^2>\boldsymbol{a}^2\boldsymbol{b}^2; \tag{13}$$

an Euclidean ***t-t*** kind or ***s-s*** kind (complex) 2-plane, if

$$D=|\langle\boldsymbol{a}|\boldsymbol{b}\rangle|^2-\boldsymbol{a}^2\boldsymbol{b}^2<0 \ \Leftrightarrow\ |\langle\boldsymbol{a}|\boldsymbol{b}\rangle|^2<\boldsymbol{a}^2\boldsymbol{b}^2; \tag{14}$$

degenerate (complex) 2-plane, if

$$D=|\langle\boldsymbol{a}|\boldsymbol{b}\rangle|^2-\boldsymbol{a}^2\boldsymbol{b}^2=0 \ \Leftrightarrow\ |\langle\boldsymbol{a}|\boldsymbol{b}\rangle|^2=\boldsymbol{a}^2\boldsymbol{b}^2. \tag{15}$$

To finalize the proof, we have to check that inverse statements are valid:

if $\boldsymbol{a},\boldsymbol{b}\in\mathbb{P}^2_{1,1}\subset\mathbb{P}^{n+m}_{n,m}$ $(\boldsymbol{a}\neq k\boldsymbol{b},\ k\in\mathbb{C})$, then (13) is correct;

if $\boldsymbol{a},\boldsymbol{b}\in\mathbb{P}^2_{2,0}\subset\mathbb{P}^{n+m}_{n,m}$ or $\boldsymbol{a},\boldsymbol{b}\in\mathbb{P}^2_{0,2}\subset\mathbb{P}^{n+m}_{n,m}$ $(\boldsymbol{a}\neq k\boldsymbol{b},\ k\in\mathbb{C})$, then (14) is correct;

if $\boldsymbol{a},\boldsymbol{b}$ $(\boldsymbol{a}\neq k\boldsymbol{b},\ k\in\mathbb{C})$ belong to any kind degenerate 2-plane, then (15) is correct.

Indeed, one has

$$\boldsymbol{a},\boldsymbol{b}\in\mathbb{P}^2_{1,1} \ \Leftrightarrow\ \boldsymbol{a}=\lambda\boldsymbol{t}+\mu\boldsymbol{s},\ \boldsymbol{b}=\rho\boldsymbol{t}+\sigma\boldsymbol{s}\ (\boldsymbol{t}^2=-\boldsymbol{s}^2=1,\ \langle\boldsymbol{t}|\boldsymbol{s}\rangle=0,\ \lambda,\mu,\rho,\sigma\in\mathbb{C}) \Rightarrow$$
$$\boldsymbol{a}^2\boldsymbol{b}^2=(\bar\lambda\lambda-\bar\mu\mu)(\bar\rho\rho-\bar\sigma\sigma),\quad |\langle\boldsymbol{a}|\boldsymbol{b}\rangle|^2=(\bar\lambda\rho-\bar\mu\sigma)(\lambda\bar\rho-\mu\bar\sigma)\quad\Rightarrow$$
$$D=|\langle\boldsymbol{a}|\boldsymbol{b}\rangle|^2-\boldsymbol{a}^2\boldsymbol{b}^2=\bar\lambda\lambda\bar\sigma\sigma+\bar\mu\mu\bar\rho\rho-\bar\lambda\mu\rho\bar\sigma-\lambda\bar\mu\bar\rho\sigma=\overline{(\lambda\sigma-\mu\rho)}(\lambda\sigma-\mu\rho)\geq 0,$$

that is (13). Besides, the equality in the last expression holds iff $\lambda/\rho=\mu/\sigma=k \Leftrightarrow \boldsymbol{a}=k\boldsymbol{b},\ k\in\mathbb{C}$.

For $\boldsymbol{a},\boldsymbol{b}\in\mathbb{P}^2_{2,0}\subset\mathbb{P}^{n+m}_{n,m}$ one has $(\boldsymbol{t}_1^2=\boldsymbol{t}_2^2=1,\ \langle\boldsymbol{t}_1|\boldsymbol{t}_2\rangle=0,\ \ \lambda,\mu,\rho,\sigma\in\mathbb{C})$

$$\boldsymbol{a}=\lambda\boldsymbol{t}_1+\mu\boldsymbol{t}_2,\ \ \boldsymbol{b}=\rho\boldsymbol{t}_1+\sigma\boldsymbol{t}_2\ \ \boldsymbol{a}^2\boldsymbol{b}^2=(\bar\lambda\lambda+\bar\mu\mu)(\bar\rho\rho+\bar\sigma\sigma),\ \ |\langle\boldsymbol{a}|\boldsymbol{b}\rangle|^2=(\bar\lambda\rho+\bar\mu\sigma)(\lambda\bar\rho+\mu\bar\sigma),$$
$$\Rightarrow D=|\langle\boldsymbol{a}|\boldsymbol{b}\rangle|^2-\boldsymbol{a}^2\boldsymbol{b}^2=\bar\lambda\rho\mu\bar\sigma+\bar\mu\sigma\lambda\bar\rho-\bar\lambda\lambda\bar\sigma\sigma-\bar\mu\mu\bar\rho\rho=-\overline{(\lambda\sigma-\mu\rho)}(\lambda\sigma-\mu\rho)\leq 0.$$

The same result one gets for $\boldsymbol{a},\boldsymbol{b}\in\mathbb{P}^2_{0,2}$ just replacing time-like vectors $\boldsymbol{t}_1,\boldsymbol{t}_2$ by space-like vectors $\boldsymbol{s}_1,\boldsymbol{s}_2$. Hence, (14) is hold, the equality being achieved iff $\lambda/\rho=\mu/\sigma=k \Leftrightarrow \boldsymbol{a}=k\boldsymbol{b},\ k\in\mathbb{C}$, just as the Theorem 2 states.

For $\boldsymbol{a},\boldsymbol{b}$ that belong to ***c-t*** kind (or to ***c-s*** kind) degenerate 2-plane,

$$\boldsymbol{a}=\lambda_1\boldsymbol{c}+\lambda_2\boldsymbol{t}(\boldsymbol{s}),\ \ \boldsymbol{b}=\mu_1\boldsymbol{c}+\mu_2\boldsymbol{t}(\boldsymbol{s})\ \ (\boldsymbol{c}^2=0,\ \langle\boldsymbol{c}|\boldsymbol{t}(\boldsymbol{s})\rangle=0,\ \lambda_i,\mu_i\in\mathbb{C},\ i=1,2) \Rightarrow$$
$$\boldsymbol{a}^2\boldsymbol{b}^2=(\pm\bar\lambda_2\lambda_2)(\pm\bar\mu_2\mu_2),\ \ |\langle\boldsymbol{a}|\boldsymbol{b}\rangle|^2=(\pm\bar\lambda_2\mu_2)(\pm\lambda_2\bar\mu_2)\ \Rightarrow\ D=|\langle\boldsymbol{a}|\boldsymbol{b}\rangle|^2-\boldsymbol{a}^2\boldsymbol{b}^2=0,$$

and for $\boldsymbol{a}\neq\boldsymbol{0},\ \boldsymbol{b}\neq\boldsymbol{0}$ that belong to ***c-c*** kind degenerate 2-plane, obviously, one has

$$\boldsymbol{a}^2=0,\ \ \boldsymbol{b}^2=0,\ \ \langle\boldsymbol{a}|\boldsymbol{b}\rangle=0\ \Rightarrow D=0,$$

so, (15) is valid as well, just as the Theorem 2 states. □

**3. Lagrange's Identity and Multi-vectors.** It is easy to check that the definition (10) for *D* can be represent (using known Lagrange's identity, see, e.g., [8], [9]) as inner square of 2-blade $\boldsymbol{a}\wedge\boldsymbol{b}$. We have

**Statement 1.** *For (complex) vectors* $\boldsymbol{a},\boldsymbol{b}\in\mathbb{P}^{n+m}_{n,m}$ *Lagrange's identity gets a form*

$$\langle\boldsymbol{a}|\boldsymbol{a}\rangle\langle\boldsymbol{b}|\boldsymbol{b}\rangle-|\langle\boldsymbol{a}|\boldsymbol{b}\rangle|^2=\langle\boldsymbol{a}\wedge\boldsymbol{b}\,|\,\boldsymbol{a}\wedge\boldsymbol{b}\rangle, \tag{16}$$

*where* $\boldsymbol{a}\wedge\boldsymbol{b}$ *denotes an exterior (wedge) product* [7], [9].

**Proof.** Setting $c_i=\bar{a}^i$ and $d_j=\bar{b}^j$ in the **Binet–Cauchy identity**

$$a_ic_ib_jd_j=a_id_ib_jc_j+\frac{1}{2}(a_ib_j-b_ia_j)(c_id_j-d_ic_j)$$

that is valid (see, e.g., [12]) for any set $a_i,c_i,b_i,d_i,\ i=\overline{1,n+m}$, of complex numbers (or, more generally, elements of a commutative ring) one gets

$$a_i\bar{a}^ib_j\bar{b}^j=a_i\bar{b}^ib_j\bar{a}^j+\frac{1}{2}(a_ib_j-b_ia_j)(\bar{a}^i\bar{b}^j-\bar{b}^i\bar{a}^j),$$

that is the Lagrange's identity for (complex) vectors $\boldsymbol{a},\boldsymbol{b}\in\mathbb{P}^{n+m}_{n,m}$. Taking here into account (4'), we get

$$g_{ik}\bar{a}^ia^kg_{jl}\bar{b}^jb^l-g_{jl}\bar{a}^jb^lg_{ik}\bar{b}^ia^k=\frac{1}{2}(\bar{a}^i\bar{b}^j-\bar{b}^i\bar{a}^j)g_{il}g_{jk}(a^lb^k-b^la^k). \tag{17}$$

According to the definition of inner product on the exterior algebra $\Lambda^k(\mathbb{V}),\ \forall k\in\mathbb{N},$ extending the inner product defined on the underlying vector space $\mathbb{V}$ (see, e. g., [7], [10])

$$\langle\wedge^k\boldsymbol{a}_i\,|\,\wedge^k\boldsymbol{b}_j\rangle=\det\langle\boldsymbol{a}_i|\boldsymbol{b}_j\rangle,\quad \boldsymbol{a}_i,\boldsymbol{b}_j\in\mathbb{V},\ i,j=\overline{1,k}, \tag{18}$$

for 2-blades – decomposable bivectors – one gets

$$\langle\boldsymbol{a}_1\wedge\boldsymbol{a}_2\,|\,\boldsymbol{b}_1\wedge\boldsymbol{b}_2\rangle=\det\begin{bmatrix}\langle\boldsymbol{a}_1|\boldsymbol{b}_1\rangle & \langle\boldsymbol{a}_2|\boldsymbol{b}_1\rangle\\ \langle\boldsymbol{a}_1|\boldsymbol{b}_2\rangle & \langle\boldsymbol{a}_2|\boldsymbol{b}_2\rangle\end{bmatrix}=\langle\boldsymbol{a}_1|\boldsymbol{b}_1\rangle\langle\boldsymbol{a}_2|\boldsymbol{b}_2\rangle-\langle\boldsymbol{a}_1|\boldsymbol{b}_2\rangle\langle\boldsymbol{a}_2|\boldsymbol{b}_1\rangle. \tag{18'}$$

Therefore, the expression (17) turns into (16) $(\boldsymbol{a}_1=\boldsymbol{b}_1=\boldsymbol{a},\ \boldsymbol{a}_2=\boldsymbol{b}_2=\boldsymbol{b}\in\mathbb{P}^{n+m}_{n,m})$. □

As far all 2-blades belong to (complex) linear vector space $\Lambda^2(\mathbb{P}^{n+m}_{n,m}),$ the formula (16) and the theorems proved above (in real and complex cases, correspondingly) implies that holds the next

**Statement 2.** *A 2-blade* $\boldsymbol{a}\wedge\boldsymbol{b}\neq\boldsymbol{0}_{(2)}$ *as a vector in* $\Lambda^2(\mathbb{P}^{n+m}_{n,m})$ *is*

*1) space-like iff* $\boldsymbol{a},\boldsymbol{b}\in\mathbb{P}^2_{1,1}$;

*2) time-like iff* $\boldsymbol{a},\boldsymbol{b}\in\mathbb{P}^2_{2,0}$ *or* $\boldsymbol{a},\boldsymbol{b}\in\mathbb{P}^2_{0,2}$;

*3) isotropic iff both vectors* $\boldsymbol{a},\boldsymbol{b}\in\mathbb{P}^{n+m}_{n,m}$ *belong to (any kind) degenerate 2-plane in* $\mathbb{P}^{n+m}_{n,m}$.

**Proof.** Due to the condition $\boldsymbol{a}\wedge\boldsymbol{b}\neq\boldsymbol{0}_{(2)}$ (zero vector in $\Lambda^2(\mathbb{P}^{n+m}_{n,m})$ hereinafter is denoted as $\boldsymbol{0}_{(k)},\ k\in\mathbb{N}$; $\boldsymbol{0}_{(1)}=\boldsymbol{0}\in\mathbb{P}^{n+m}_{n,m}$) we have $\boldsymbol{a}\neq\alpha\boldsymbol{b},\ \alpha\in\mathbb{C}$. Therefore, if $\boldsymbol{a},\boldsymbol{b}\in\mathbb{P}^2_{1,1}$ then, according to (16) and (13), one has

$$\langle\boldsymbol{a}\wedge\boldsymbol{b}\,|\,\boldsymbol{a}\wedge\boldsymbol{b}\rangle=\langle\boldsymbol{a}|\boldsymbol{a}\rangle\langle\boldsymbol{b}|\boldsymbol{b}\rangle-|\langle\boldsymbol{a}|\boldsymbol{b}\rangle|^2<0,$$

so, the condition $\boldsymbol{a},\boldsymbol{b}\in\mathbb{P}^2_{1,1}$ in 1) is sufficient;

if $\boldsymbol{a},\boldsymbol{b}\in\mathbb{P}^2_{2,0}$ (or if $\boldsymbol{a},\boldsymbol{b}\in\mathbb{P}^2_{0,2}$) then (16) and (14) give

$$\langle\boldsymbol{a}\wedge\boldsymbol{b}\,|\,\boldsymbol{a}\wedge\boldsymbol{b}\rangle=\langle\boldsymbol{a}|\boldsymbol{a}\rangle\langle\boldsymbol{b}|\boldsymbol{b}\rangle-|\langle\boldsymbol{a}|\boldsymbol{b}\rangle|^2>0$$

and hence the condition $\boldsymbol{a},\boldsymbol{b}\in\mathbb{P}^2_{2,0}$ (or $\boldsymbol{a},\boldsymbol{b}\in\mathbb{P}^2_{0,2}$) in 2) is sufficient as well.

Finally, if both vectors $\boldsymbol{a}, \boldsymbol{b} \in \mathbb{P}^{n+m}_{n,m}$ belong to degenerate 2-plane (in $\mathbb{P}^{n+m}_{n,m}$) of any kind, then, according to (16) and (15), one has

$$\langle \boldsymbol{a}\wedge\boldsymbol{b} | \boldsymbol{a}\wedge\boldsymbol{b}\rangle = \langle \boldsymbol{a}|\boldsymbol{a}\rangle\langle \boldsymbol{b}|\boldsymbol{b}\rangle - |\langle \boldsymbol{a}|\boldsymbol{b}\rangle|^2 = 0,$$

so, the condition above is sufficient as well.

Clearly, according to the Theorem 2, the Statement's conditions 1)-3) are necessary, as well. □

**Corollary 1.** The linear (complex) vector space $\Lambda^2(\mathbb{P}^{n+m}_{n,m})$ is a pseudo-Euclidean $\mathbb{P}^{N+M}_{N,M}$ space with Cartesian metric

$$G_{ij} = G^{ij} = \mathrm{diag}(\overbrace{+1,...,+1}^{N}, \overbrace{-1,...,-1}^{M}), \quad i,j = \overline{1, N+M}, \tag{19}$$

and its dimensions are

$$N = \binom{n}{2} + \binom{m}{2}, \; M = nm \;\Rightarrow\; N+M = \binom{m+n}{2}. \tag{20}$$

**Proof.** According to Statement 2, basic 2-blades $\boldsymbol{E}_i = \boldsymbol{e}_{i_1}\wedge\boldsymbol{e}_{i_2}$ and $\boldsymbol{E}_j = \boldsymbol{e}_{j_1}\wedge\boldsymbol{e}_{j_2}$ are orthonormal with respect to the inner product of 2-blades (18') and relations (4'):

$$G_{ij} = \langle \boldsymbol{E}_i | \boldsymbol{E}_j \rangle = \langle \boldsymbol{e}_{i_1}\wedge\boldsymbol{e}_{i_2} | \boldsymbol{e}_{j_1}\wedge\boldsymbol{e}_{j_2} \rangle = \langle \boldsymbol{e}_{i_1}|\boldsymbol{e}_{j_1}\rangle\langle \boldsymbol{e}_{i_2}|\boldsymbol{e}_{j_2}\rangle - \langle \boldsymbol{e}_{i_1}|\boldsymbol{e}_{j_2}\rangle\langle \boldsymbol{e}_{i_2}|\boldsymbol{e}_{j_1}\rangle = \det\left[g_{i_\alpha j_\beta}\right]_{\alpha,\beta=1,2}, \tag{21}$$

and satisfy (19) for corresponding numeration of 2-blades $\boldsymbol{E}_i$. Then, simple combinatorics gives us (20).

**Corollary 2.** In pseudo-Euclidean (complex) space $\Lambda^2(\mathbb{P}^{n+m}_{n,m})$ the Theorem 2 is valid, that gives the next inequalities ($\boldsymbol{A}, \boldsymbol{B}$ being any vectors in $\Lambda^2(\mathbb{P}^{n+m}_{n,m})$, not necessary 2-blades, $\boldsymbol{A} \neq k\boldsymbol{B}$, $k \in \mathbb{R}$):

$$\langle \boldsymbol{A}|\boldsymbol{A}\rangle\langle \boldsymbol{B}|\boldsymbol{B}\rangle - |\langle \boldsymbol{A}|\boldsymbol{B}\rangle|^2 < 0, \tag{22$_1$}$$

iff span$[\boldsymbol{A},\boldsymbol{B}]$ is pseudo-Euclidean "bi-plane";

$$\langle \boldsymbol{A}|\boldsymbol{A}\rangle\langle \boldsymbol{B}|\boldsymbol{B}\rangle - |\langle \boldsymbol{A}|\boldsymbol{B}\rangle|^2 > 0, \tag{22$_2$}$$

iff span$[\boldsymbol{A},\boldsymbol{B}]$ is Euclidean "bi-plane";

$$\langle \boldsymbol{A}|\boldsymbol{A}\rangle\langle \boldsymbol{B}|\boldsymbol{B}\rangle - |\langle \boldsymbol{A}|\boldsymbol{B}\rangle|^2 = 0, \tag{22$_3$}$$

iff span$[\boldsymbol{A},\boldsymbol{B}]$ is isotropic "bi-plane".

Applying the formulas (22$_{1\text{-}3}$) to the particular case of 2-blades $\boldsymbol{A} = \boldsymbol{a}\wedge\boldsymbol{b}$, $\boldsymbol{B} = \boldsymbol{c}\wedge\boldsymbol{d}$ we get criteria for span$[\boldsymbol{a}\wedge\boldsymbol{b}, \boldsymbol{c}\wedge\boldsymbol{d}]$ to be pseudo-Euclidean, Euclidean or isotropic "bi-planes", respectively:

$$\det\begin{bmatrix} \boldsymbol{a}^2\boldsymbol{b}^2 - |\langle \boldsymbol{a}|\boldsymbol{b}\rangle|^2 & \langle \boldsymbol{a}|\boldsymbol{c}\rangle\langle \boldsymbol{b}|\boldsymbol{d}\rangle - \langle \boldsymbol{a}|\boldsymbol{d}\rangle\langle \boldsymbol{b}|\boldsymbol{c}\rangle \\ \langle \boldsymbol{c}|\boldsymbol{a}\rangle\langle \boldsymbol{d}|\boldsymbol{b}\rangle - \langle \boldsymbol{d}|\boldsymbol{a}\rangle\langle \boldsymbol{c}|\boldsymbol{b}\rangle & \boldsymbol{c}^2\boldsymbol{d}^2 - |\langle \boldsymbol{c}|\boldsymbol{d}\rangle|^2 \end{bmatrix} < 0, \tag{23$_1$}$$

$$\det\begin{bmatrix} \boldsymbol{a}^2\boldsymbol{b}^2 - |\langle \boldsymbol{a}|\boldsymbol{b}\rangle|^2 & \langle \boldsymbol{a}|\boldsymbol{c}\rangle\langle \boldsymbol{b}|\boldsymbol{d}\rangle - \langle \boldsymbol{a}|\boldsymbol{d}\rangle\langle \boldsymbol{b}|\boldsymbol{c}\rangle \\ \langle \boldsymbol{c}|\boldsymbol{a}\rangle\langle \boldsymbol{d}|\boldsymbol{b}\rangle - \langle \boldsymbol{d}|\boldsymbol{a}\rangle\langle \boldsymbol{c}|\boldsymbol{b}\rangle & \boldsymbol{c}^2\boldsymbol{d}^2 - |\langle \boldsymbol{c}|\boldsymbol{d}\rangle|^2 \end{bmatrix} > 0, \tag{23$_2$}$$

$$\det\begin{bmatrix} \boldsymbol{a}^2\boldsymbol{b}^2 - |\langle \boldsymbol{a}|\boldsymbol{b}\rangle|^2 & \langle \boldsymbol{a}|\boldsymbol{c}\rangle\langle \boldsymbol{b}|\boldsymbol{d}\rangle - \langle \boldsymbol{a}|\boldsymbol{d}\rangle\langle \boldsymbol{b}|\boldsymbol{c}\rangle \\ \langle \boldsymbol{c}|\boldsymbol{a}\rangle\langle \boldsymbol{d}|\boldsymbol{b}\rangle - \langle \boldsymbol{d}|\boldsymbol{a}\rangle\langle \boldsymbol{c}|\boldsymbol{b}\rangle & \boldsymbol{c}^2\boldsymbol{d}^2 - |\langle \boldsymbol{c}|\boldsymbol{d}\rangle|^2 \end{bmatrix} = 0. \tag{23$_3$}$$

**4. Angles and some generalizations.** Proved Theorem 1 allows one to generalize the concept of an angle between vectors $\boldsymbol{a},\boldsymbol{b}\in\mathbb{P}^{n+m}_{n,m}$. For a case of real space $\mathbb{P}^{n+m}_{n,m}$ one has:

*a*) If vectors $\boldsymbol{a},\boldsymbol{b}\in\mathbb{P}^{2}_{2,0}=\mathbb{E}^2\subset\mathbb{P}^{n+m}_{n,m}$ or $\boldsymbol{a},\boldsymbol{b}\in\mathbb{P}^{2}_{0,2}=\mathbb{E}^2\subset\mathbb{P}^{n+m}_{n,m}$, then $\exists\varphi_{(\boldsymbol{a},\boldsymbol{b})}\in\mathbb{R}$, such that

$$\cos\varphi_{(\boldsymbol{a},\boldsymbol{b})}=-\langle\boldsymbol{a}|\boldsymbol{b}\rangle(\boldsymbol{a}^2\boldsymbol{b}^2)^{-1/2},\ \boldsymbol{a}^2,\boldsymbol{b}^2<0;\quad \cos\varphi_{(\boldsymbol{a},\boldsymbol{b})}=\langle\boldsymbol{a}|\boldsymbol{b}\rangle(\boldsymbol{a}^2\boldsymbol{b}^2)^{-1/2},\ \boldsymbol{a}^2,\boldsymbol{b}^2>0. \quad (24)$$

The signs that are chosen in (24) provide suitable result

$$\boldsymbol{b}=k\boldsymbol{a}\Rightarrow\cos\varphi_{(\boldsymbol{a},\boldsymbol{b})}=k/|k|,\ k\in\mathbb{R}$$

(see, e. g., [6], [4], [5]).

*b*) According to the Theorem 1, $\langle\boldsymbol{a}|\boldsymbol{b}\rangle^2>\boldsymbol{a}^2\boldsymbol{b}^2\Leftrightarrow\boldsymbol{a},\boldsymbol{b}\in\mathbb{P}^2_{1,1}\subset\mathbb{P}^{n+m}_{n,m}$, $\boldsymbol{a}\neq k\boldsymbol{b}$ then $(\boldsymbol{t}^2=-\boldsymbol{s}^2=1,\langle\boldsymbol{t}|\boldsymbol{s}\rangle=0)$

$$\boldsymbol{a}^2>0,\ \boldsymbol{b}^2<0\ \Leftrightarrow\ \boldsymbol{a}_\pm=|\boldsymbol{a}|(\pm\boldsymbol{t}\cosh\alpha+\boldsymbol{s}\sinh\alpha),\quad \boldsymbol{b}_\pm=|\boldsymbol{b}|(\boldsymbol{t}\sinh\beta\pm\boldsymbol{s}\cosh\beta)\quad\Rightarrow$$
$$\langle\boldsymbol{a}_+|\boldsymbol{b}_\pm\rangle=|\boldsymbol{a}||\boldsymbol{b}|\sinh(\beta\mp\alpha),\quad \langle\boldsymbol{a}_-|\boldsymbol{b}_\pm\rangle=-|\boldsymbol{a}||\boldsymbol{b}|\sinh(\beta\pm\alpha);$$
$$\boldsymbol{a}^2>0,\ \boldsymbol{b}^2>0\ \Leftrightarrow\ \boldsymbol{a}_\pm=|\boldsymbol{a}|(\pm\boldsymbol{t}\cosh\alpha+\boldsymbol{s}\sinh\alpha),\quad \boldsymbol{b}_\pm=|\boldsymbol{b}|(\pm\boldsymbol{t}\cosh\beta+\boldsymbol{s}\sinh\beta)\quad\Rightarrow$$
$$\langle\boldsymbol{a}_+|\boldsymbol{b}_\pm\rangle=\pm|\boldsymbol{a}||\boldsymbol{b}|\cosh(\alpha\mp\beta),\quad \langle\boldsymbol{a}_-|\boldsymbol{b}_\pm\rangle=\mp|\boldsymbol{a}||\boldsymbol{b}|\cosh(\alpha\pm\beta);$$
$$\boldsymbol{a}^2<0,\ \boldsymbol{b}^2<0\ \Leftrightarrow\ \boldsymbol{a}_\pm=|\boldsymbol{a}|(\boldsymbol{t}\sinh\alpha\pm\boldsymbol{s}\cosh\alpha),\quad \boldsymbol{b}_\pm=|\boldsymbol{b}|(\boldsymbol{t}\sinh\beta\pm\boldsymbol{s}\cosh\beta)\quad\Rightarrow$$
$$\langle\boldsymbol{a}_+|\boldsymbol{b}_\pm\rangle=\mp|\boldsymbol{a}||\boldsymbol{b}|\cosh(\alpha\mp\beta),\quad \langle\boldsymbol{a}_-|\boldsymbol{b}_\pm\rangle=\pm|\boldsymbol{a}||\boldsymbol{b}|\cosh(\alpha\pm\beta).$$
$$(\forall\alpha,\beta\in\mathbb{R})$$

If one of the vectors $\boldsymbol{a},\boldsymbol{b}\in\mathbb{P}^2_{1,1}$ is time-like while another one is space-like (hence, if $\boldsymbol{a}^2\boldsymbol{b}^2<0$), then

$$\langle\boldsymbol{a}_+|\boldsymbol{b}_\pm\rangle=|\boldsymbol{a}||\boldsymbol{b}|\sinh(\beta\mp\alpha)=0\Rightarrow\beta=\pm\alpha,\ \langle\boldsymbol{a}_-|\boldsymbol{b}_\pm\rangle=-|\boldsymbol{a}||\boldsymbol{b}|\sinh(\beta\pm\alpha)=0\ \Rightarrow\ \beta=\mp\alpha.$$

Therefore, in this case one can introduce a "defect of orthogonality" $\delta_{(\boldsymbol{a},\boldsymbol{b})}\in\mathbb{R}$ as (see also [4], [5])

$$\sinh\delta_{(\boldsymbol{a},\boldsymbol{b})}=\langle\boldsymbol{a}|\boldsymbol{b}\rangle(-\boldsymbol{a}^2\boldsymbol{b}^2)^{-1/2}, \quad (25_1)$$

Note, that determined "defect of orthogonality" δ is related with the longitudinal components of the time-like and the space-like vectors $\boldsymbol{a},\boldsymbol{b}\in\mathbb{P}^2_{1,1}\subset\mathbb{P}^{n+m}_{n,m}$: $(\boldsymbol{a}^2\boldsymbol{b}^2<0)$

$$\boldsymbol{a}_{\|}^{(\boldsymbol{b})}\equiv\frac{\langle\boldsymbol{a}|\boldsymbol{b}\rangle\boldsymbol{b}}{\boldsymbol{b}^2},\quad \boldsymbol{b}_{\|}^{(\boldsymbol{a})}\equiv\frac{\langle\boldsymbol{a}|\boldsymbol{b}\rangle\boldsymbol{a}}{\boldsymbol{a}^2}\quad\Rightarrow\quad\langle\boldsymbol{a}_{\|}^{(\boldsymbol{b})}|\boldsymbol{b}_{\|}^{(\boldsymbol{a})}\rangle=\frac{\langle\boldsymbol{a}|\boldsymbol{b}\rangle^3}{\boldsymbol{a}^2\boldsymbol{b}^2}=-|\boldsymbol{a}||\boldsymbol{b}|(\sinh\delta_{(\boldsymbol{a},\boldsymbol{b})})^3.$$

If the vectors $\boldsymbol{a},\boldsymbol{b}\in\mathbb{P}^2_{1,1}\subset\mathbb{P}^{n+m}_{n,m}$ both are ***t***-like or both are ***s***-like $(\boldsymbol{a}^2\boldsymbol{b}^2>0)$, then the "angle" $\varphi_{(\boldsymbol{a},\boldsymbol{b})}=|\alpha-\beta|\in\mathbb{R}$ between them indicates a "deviation of collinearity" of these vectors (see also [4], [5]):

$$\cosh\varphi_{(\boldsymbol{a},\boldsymbol{b})}=|\langle\boldsymbol{a}|\boldsymbol{b}\rangle|(\boldsymbol{a}^2\boldsymbol{b}^2)^{-1/2}\geq1. \quad (25_2)$$

An angle between isotropic vector and time-like (space-like) vector is not defined.

*c*) If $\langle\boldsymbol{a}|\boldsymbol{b}\rangle^2=\boldsymbol{a}^2\boldsymbol{b}^2>0$, $\boldsymbol{a}\neq k\boldsymbol{b}$ then the span[***a***, ***b***] is degenerate 2-plane, and for a ***c***-***t*** 2-plane one has

$$\boldsymbol{a}=\alpha\boldsymbol{c}+\beta\boldsymbol{t},\quad \boldsymbol{b}=\lambda\boldsymbol{c}+\mu\boldsymbol{t}\ \Rightarrow\ \cos\varphi_{(\boldsymbol{a},\boldsymbol{b})}=\langle\boldsymbol{a}|\boldsymbol{b}\rangle(\boldsymbol{a}^2\boldsymbol{b}^2)^{-1/2}=\pm1\quad(\forall\alpha,\beta,\lambda,\mu\in\mathbb{R}). \quad (26_1)$$

Correspondingly, for a ***c***-***s*** 2-plane one has

$$\boldsymbol{a}=\alpha\boldsymbol{c}+\beta\boldsymbol{s},\ \boldsymbol{b}=\lambda\boldsymbol{c}+\mu\boldsymbol{s}\ \Rightarrow\ \cos\varphi_{(\boldsymbol{a},\boldsymbol{b})}=-\langle\boldsymbol{a}|\boldsymbol{b}\rangle(\boldsymbol{a}^2\boldsymbol{b}^2)^{-1/2}=\pm1\quad(\forall\alpha,\beta,\lambda,\mu\in\mathbb{R}). \quad (26_2)$$

Obviously, these definitions give (including the limit cases $\boldsymbol{a}=\boldsymbol{b}\to\boldsymbol{c}$, see Appendix 2)

$$\cos\varphi=1\ \text{for}\ \boldsymbol{a}=\boldsymbol{b},\quad \cos\varphi=-1\ \text{for}\ \boldsymbol{a}=-\boldsymbol{b}.$$

If span[***a***, ***b***] is a ***c***-***c*** degenerate 2-plane, then an angle between vectors belong to it is not defined.

**Note**. The deviation of collinearity ($25_2$), as well as the defect of orthogonality ($25_1$) tends to infinity, as one of the vectors $\boldsymbol{a},\boldsymbol{b}\in\mathbb{P}^2_{1,1}$ tends to isotropic vector $\boldsymbol{c}\in\text{span}[\boldsymbol{a},\boldsymbol{b}]=\mathbb{P}^2_{1,1}\subset\mathbb{P}^{n+m}_{n,m}$, $\boldsymbol{c}^2=0, \boldsymbol{c}\neq\boldsymbol{0}$.

**Corollary 3.** Obviously, formulas (24)-($26_2$) are applicable to the bivectors, considering as vectors in linear real space $\Lambda^2(\mathbb{P}^{n+m}_{n,m})$. Thus, we get $(\varphi,\delta\in\mathbb{R})$:

$$\cos\varphi=-\langle\boldsymbol{A}|\boldsymbol{B}\rangle(\boldsymbol{A}^2\boldsymbol{B}^2)^{-1/2},\quad \boldsymbol{A}^2,\boldsymbol{B}^2<0;\quad \cos\varphi=\langle\boldsymbol{A}|\boldsymbol{B}\rangle(\boldsymbol{A}^2\boldsymbol{B}^2)^{-1/2},\quad \boldsymbol{A}^2,\boldsymbol{B}^2>0, \tag{27}$$

if real bivectors $\boldsymbol{A},\boldsymbol{B}$ spanned Euclidean "bi-plane";

$$\delta=\sinh^{-1}\langle\boldsymbol{A}|\boldsymbol{B}\rangle(-\boldsymbol{A}^2\boldsymbol{B}^2)^{-1/2} \tag{$28_1$}$$

if real bivectors $\boldsymbol{A},\boldsymbol{B}$ spanned pseudo-Euclidean "bi-plane" and $\boldsymbol{A}^2\boldsymbol{B}^2<0$;

$$\cosh\varphi=|\langle\boldsymbol{A}|\boldsymbol{B}\rangle|(\boldsymbol{A}^2\boldsymbol{B}^2)^{-1/2}\geq1, \tag{$28_2$}$$

if real bivectors $\boldsymbol{A},\boldsymbol{B}$ spanned pseudo-Euclidean "bi-plane" and $\boldsymbol{A}^2\boldsymbol{B}^2>0$;

$$\cos\varphi=\langle\boldsymbol{A}|\boldsymbol{B}\rangle(\boldsymbol{A}^2\boldsymbol{B}^2)^{-1/2}=\pm1, \tag{$29_1$}$$

if real bivectors $\boldsymbol{A},\boldsymbol{B}$ spanned ***c-t*** kind degenerate "bi-plane" and

$$\cos\varphi=-\langle\boldsymbol{A}|\boldsymbol{B}\rangle(\boldsymbol{A}^2\boldsymbol{B}^2)^{-1/2}=\pm1, \tag{$29_2$}$$

if real bivectors $\boldsymbol{A},\boldsymbol{B}$ spanned ***c-s*** kind degenerate "bi-plane".

In particular, specifying (real) bivectors $\boldsymbol{A},\boldsymbol{B}$ as 2-blades $\boldsymbol{A}=\boldsymbol{a}_1\wedge\boldsymbol{a}_2$, $\boldsymbol{B}=\boldsymbol{a}_3\wedge\boldsymbol{a}_4$, $\boldsymbol{A}\neq k\boldsymbol{B}$ $(\forall k\in\mathbb{R})$, and calculating inner products according (18') one gets some generalizations of dihedral angles [11]. To show this, we need previously the next Lemmas:

**Lemma 1.** *If* $\boldsymbol{a}_1\wedge\boldsymbol{a}_2$, $\boldsymbol{a}_3\wedge\boldsymbol{a}_4\in\Lambda^2\mathbb{V}$ *are nonzero 2-blades,* $\boldsymbol{a}_1\wedge\boldsymbol{a}_2\neq k(\boldsymbol{a}_3\wedge\boldsymbol{a}_4)$, $k\in\mathbb{R}$, $\boldsymbol{a}_j$, $j=\overline{1,4}$, *being vectors in a linear real space* $\mathbb{V}(\mathbb{R})$, *then the next statement is valid:*

$$\boldsymbol{a}_1\wedge\boldsymbol{a}_2\wedge\boldsymbol{a}_3\wedge\boldsymbol{a}_4=\boldsymbol{0}_{(4)}\ \Leftrightarrow\ \exists\alpha,\beta,\lambda,\mu\in\mathbb{R},\ \boldsymbol{x}=\alpha\boldsymbol{a}_1+\beta\boldsymbol{a}_2=\lambda\boldsymbol{a}_3+\mu\boldsymbol{a}_4\neq\boldsymbol{0}. \tag{30}$$

**Proof.** *a*) Geometric approach: the condition $\boldsymbol{a}_1\wedge\boldsymbol{a}_2\wedge\boldsymbol{a}_3\wedge\boldsymbol{a}_4=\Lambda^4\boldsymbol{a}_j=\boldsymbol{0}_{(4)}$ means that 2-blades $\boldsymbol{a}_1\wedge\boldsymbol{a}_2$, $\boldsymbol{a}_3\wedge\boldsymbol{a}_4$ belong to at most 3-dimensional sub-space of *n*-dimensional linear space $\mathbb{V}(\mathbb{R})$ while the condition $\boldsymbol{a}_1\wedge\boldsymbol{a}_2\neq k(\boldsymbol{a}_3\wedge\boldsymbol{a}_4)$, $k\in\mathbb{R}$, indicates that the dimension of this sub-space is at least 3. Thus, both 2-blades $\boldsymbol{a}_1\wedge\boldsymbol{a}_2$, $\boldsymbol{a}_3\wedge\boldsymbol{a}_4$ belong to the same 3-dimensional linear space. Therefore, non-parallel 2-planes $\text{span}[\boldsymbol{a}_1,\boldsymbol{a}_2]$ and $\text{span}[\boldsymbol{a}_3,\boldsymbol{a}_4]$ intersect along a vector $\boldsymbol{x}$,

$$\boldsymbol{x}=\alpha\boldsymbol{a}_1+\beta\boldsymbol{a}_2=\lambda\boldsymbol{a}_3+\mu\boldsymbol{a}_4\neq\boldsymbol{0},\ \alpha,\beta,\lambda,\mu\in\mathbb{R}.$$

So, (30) is valid.

*b*) Algebraic approach: the components of the 4-blade $\Lambda^4\boldsymbol{a}_k=\boldsymbol{0}_{(4)}$, $k=\overline{1,4}$, are minors of $4\times n$ matrix $[a_k^{\ i_k}]^{k=\overline{1,4}}_{i_k=\overline{1,n}}$. So, according to the conditions of the Lemma, for the rank *r* of this matrix one has $2\leq r\leq3$. But, if $r=2$, then all 3-minors must be zero and thus all 2-minors must be proportional to each other, in contradiction with the condition $\boldsymbol{a}_1\wedge\boldsymbol{a}_2\neq k(\boldsymbol{a}_3\wedge\boldsymbol{a}_4)$, $k\in\mathbb{R}$. Hence, $r=3$ and the vectors $\boldsymbol{a}_k$, $k=\overline{1,4}$, are linearly dependent: there exist real constants $\alpha,\beta,\lambda,\mu$, not all simultaneously being zero, such that

$$\alpha\boldsymbol{a}_1+\beta\boldsymbol{a}_2+\lambda\boldsymbol{a}_3+\mu\boldsymbol{a}_4=\boldsymbol{0}\ \Leftrightarrow\ \alpha a_1^{\ i}+\beta a_2^{\ i}+\lambda a_3^{\ i}+\mu a_4^{\ i}=0,\ \ i=\overline{1,n}. \tag{31}$$

The system of equations (31) (with rank $r=3$) has the solution $\tilde\alpha,\tilde\beta,\tilde\lambda,\tilde\mu$, (unique up to proportionality) and therefore there does exist the vector $\boldsymbol{x}\neq\boldsymbol{0}$ that is common for the $\text{span}[\boldsymbol{a}_1,\boldsymbol{a}_2]$ and $\text{span}[\boldsymbol{a}_3,\boldsymbol{a}_4]$:

$$\text{span}[\boldsymbol{a}_1,\boldsymbol{a}_2]\cap\text{span}[\boldsymbol{a}_3,\boldsymbol{a}_4]=\boldsymbol{x}=\tilde\alpha\boldsymbol{a}_1+\tilde\beta\boldsymbol{a}_2=-\tilde\lambda\boldsymbol{a}_3-\tilde\mu\boldsymbol{a}_4\neq\boldsymbol{0}. \tag{32}$$

Obviously, the condition (32) is sufficient as well for the equality $\Lambda^4 \boldsymbol{a}_k = \boldsymbol{0}_{(4)},\ k=\overline{1,4}$:

$$\tilde{\alpha}\boldsymbol{a}_1 + \tilde{\beta}\boldsymbol{a}_2 = -\tilde{\lambda}\boldsymbol{a}_3 - \tilde{\mu}\boldsymbol{a}_4 \neq \boldsymbol{0} \Rightarrow \begin{cases} \tilde{\alpha}\boldsymbol{a}_1 \wedge (\boldsymbol{a}_2 \wedge \boldsymbol{a}_3 \wedge \boldsymbol{a}_4) = \tilde{\alpha}(\wedge^4 \boldsymbol{a}_k) = \boldsymbol{0}_{(4)}, \\ \tilde{\beta}\boldsymbol{a}_2 \wedge (\boldsymbol{a}_1 \wedge \boldsymbol{a}_3 \wedge \boldsymbol{a}_4) = -\tilde{\beta}(\wedge^4 \boldsymbol{a}_k) = \boldsymbol{0}_{(4)}, \end{cases} \quad \tilde{\alpha}^2 + \tilde{\beta}^2 \neq 0$$

$$\Rightarrow \boldsymbol{a}_1 \wedge \boldsymbol{a}_2 \wedge \boldsymbol{a}_3 \wedge \boldsymbol{a}_4 = \boldsymbol{0}_{(4)}. \qquad \square$$

**Corollary 4.** From (32), without loss of generality, one can find (taking into account (18))

$$\begin{aligned} &\tilde{\beta}\boldsymbol{a}_1 \wedge \boldsymbol{a}_2 \wedge \boldsymbol{a}_3 = \tilde{\mu}\boldsymbol{a}_1 \wedge \boldsymbol{a}_3 \wedge \boldsymbol{a}_4, \quad \tilde{\alpha}\boldsymbol{a}_1 \wedge \boldsymbol{a}_2 \wedge \boldsymbol{a}_4 = \tilde{\lambda}\boldsymbol{a}_2 \wedge \boldsymbol{a}_3 \wedge \boldsymbol{a}_4, \\ &\tilde{\alpha}\boldsymbol{a}_1 \wedge \boldsymbol{a}_3 \wedge \boldsymbol{a}_4 = -\tilde{\beta}\boldsymbol{a}_2 \wedge \boldsymbol{a}_3 \wedge \boldsymbol{a}_4, \quad \tilde{\lambda}\boldsymbol{a}_1 \wedge \boldsymbol{a}_2 \wedge \boldsymbol{a}_3 = -\tilde{\mu}\boldsymbol{a}_1 \wedge \boldsymbol{a}_2 \wedge \boldsymbol{a}_4. \end{aligned} \Rightarrow$$

$$\frac{\tilde{\alpha}^2}{\tilde{\lambda}^2} = \frac{(\boldsymbol{a}_2 \wedge \boldsymbol{a}_3 \wedge \boldsymbol{a}_4)^2}{(\boldsymbol{a}_1 \wedge \boldsymbol{a}_2 \wedge \boldsymbol{a}_4)^2}, \quad \frac{\tilde{\beta}^2}{\tilde{\lambda}^2} = \frac{(\boldsymbol{a}_1 \wedge \boldsymbol{a}_3 \wedge \boldsymbol{a}_4)^2}{(\boldsymbol{a}_1 \wedge \boldsymbol{a}_2 \wedge \boldsymbol{a}_4)^2}, \quad \frac{\tilde{\mu}^2}{\tilde{\lambda}^2} = \frac{(\boldsymbol{a}_1 \wedge \boldsymbol{a}_2 \wedge \boldsymbol{a}_3)^2}{(\boldsymbol{a}_1 \wedge \boldsymbol{a}_2 \wedge \boldsymbol{a}_4)^2} \Rightarrow$$

$$\begin{aligned} \boldsymbol{x} &= \sqrt{(\boldsymbol{a}_2 \wedge \boldsymbol{a}_3 \wedge \boldsymbol{a}_4)^2}\,\boldsymbol{a}_1 - \sqrt{(\boldsymbol{a}_1 \wedge \boldsymbol{a}_3 \wedge \boldsymbol{a}_4)^2}\,\boldsymbol{a}_2 = \\ &= -\sqrt{(\boldsymbol{a}_1 \wedge \boldsymbol{a}_2 \wedge \boldsymbol{a}_4)^2}\,\boldsymbol{a}_3 + \sqrt{(\boldsymbol{a}_1 \wedge \boldsymbol{a}_2 \wedge \boldsymbol{a}_3)^2}\,\boldsymbol{a}_4. \end{aligned} \tag{33}$$

**Corollary 5.** According to the Lemma 1,

$$\text{span}[\boldsymbol{a}_1, \boldsymbol{a}_2] \cap \text{span}[\boldsymbol{a}_3, \boldsymbol{a}_4] \neq \boldsymbol{0} \Leftrightarrow \Lambda^4 \boldsymbol{a}_k = \boldsymbol{0}_{(4)}, k=\overline{1,4},$$

hence

$$\boldsymbol{a}_1 \wedge \boldsymbol{a}_2 \wedge \boldsymbol{a}_3 \wedge \boldsymbol{a}_4 \neq \boldsymbol{0}_{(4)} \Leftrightarrow \text{span}[\boldsymbol{a}_1, \boldsymbol{a}_2] \cap \text{span}[\boldsymbol{a}_3, \boldsymbol{a}_4] = \boldsymbol{0} \Rightarrow$$

$$\langle \boldsymbol{a}_1 \wedge \boldsymbol{a}_2 | \boldsymbol{a}_3 \wedge \boldsymbol{a}_4 \rangle = 0. \tag{34}$$

**Lemma 2.** *In a pseudo-Euclidean* $\mathbb{P}^{n+m}_{n,m}$ *(non-degenerate) space a normal vector* $\boldsymbol{n}^{(k)} \in \mathbb{P}^{n+m}_{n,m}$ *to* $n+m-1$-*dimensional hyper-plane* $\text{span}[\boldsymbol{a}_1, \boldsymbol{a}_2, ..., \not{\boldsymbol{a}}_k, ..., \boldsymbol{a}_{n+m}] \subset \mathbb{P}^{n+m}_{n,m}$ *can be presented as*

$$\boldsymbol{n}^{(k)} = \alpha^{(k)} \frac{\det\left[A^{(k)}\right]}{\det[A]} = \alpha^{(k)} \frac{\boldsymbol{a}_l \Delta^{kl}}{\det[A]}, \ \alpha^{(k)} \in \mathbb{R}, \quad k = \overline{1, n+m}, \tag{35}$$

*where* $[A]$ *is Gram's matrix formed by scalar products* $\langle \boldsymbol{a}_j | \boldsymbol{a}_k \rangle,\ j,k = \overline{1,n+m}$, *and* $\left[A^{(k)}\right]$ *– the matrix obtained from it with replacement of k-th row by (linearly independent) vectors* $\boldsymbol{a}_1, \boldsymbol{a}_2, ..., \boldsymbol{a}_{n+m} \in \mathbb{P}^{n+m}_{n,m}$:

$$[A] = \left[\langle \boldsymbol{a}_j | \boldsymbol{a}_k \rangle\right]_{j,k=1,n+m}, \quad \left[A^{(k)}\right] = \begin{bmatrix} \langle \boldsymbol{a}_1 | \boldsymbol{a}_1 \rangle & \langle \boldsymbol{a}_1 | \boldsymbol{a}_2 \rangle & \cdots & \langle \boldsymbol{a}_1 | \boldsymbol{a}_{n+m} \rangle \\ \langle \boldsymbol{a}_2 | \boldsymbol{a}_1 \rangle & \langle \boldsymbol{a}_2 | \boldsymbol{a}_2 \rangle & \cdots & \langle \boldsymbol{a}_2 | \boldsymbol{a}_{n+m} \rangle \\ \vdots & \vdots & \cdots & \vdots \\ \langle \boldsymbol{a}_{k-1} | \boldsymbol{a}_1 \rangle & \langle \boldsymbol{a}_{k-1} | \boldsymbol{a}_2 \rangle & \cdots & \langle \boldsymbol{a}_{k-1} | \boldsymbol{a}_{n+m} \rangle \\ \boldsymbol{a}_1 & \boldsymbol{a}_2 & \cdots & \boldsymbol{a}_{n+m} \\ \langle \boldsymbol{a}_{k+1} | \boldsymbol{a}_1 \rangle & \langle \boldsymbol{a}_{k+1} | \boldsymbol{a}_2 \rangle & \cdots & \langle \boldsymbol{a}_{k+1} | \boldsymbol{a}_{n+m} \rangle \\ \vdots & \vdots & \cdots & \vdots \\ \langle \boldsymbol{a}_{n+m} | \boldsymbol{a}_1 \rangle & \langle \boldsymbol{a}_{n+m} | \boldsymbol{a}_2 \rangle & \cdots & \langle \boldsymbol{a}_{n+m} | \boldsymbol{a}_{n+m} \rangle \end{bmatrix}, \tag{36}$$

*and* $\Delta^{kl}$ *denotes the algebraic complement to* [A] *matrix's entry* $\langle \boldsymbol{a}_k | \boldsymbol{a}_l \rangle = A_{kl},\ k,l = \overline{1,n+m}$.

**Proof.** Since Gram's matrix $[A]$ is non-singular, $A = \det[A] \neq 0$, then direct calculation shows that

$$\langle \boldsymbol{a}_j | \boldsymbol{n}^{(k)} \rangle = \alpha^{(k)} \delta^k_j, \quad j,k = \overline{1,n+m}. \qquad \square$$

**Corollary 3.** Based on Laplace Theorem it is easy to calculate

$$\left\langle \boldsymbol{n}^{(k)} \middle| \boldsymbol{n}^{(j)} \right\rangle = \alpha^{(k)}\alpha^{(j)} \frac{\left\langle \det\left[A^{(k)}\right] \middle| \det\left[A^{(j)}\right] \right\rangle}{\det^2[A]} = \alpha^{(k)}\alpha^{(j)} \frac{\left\langle \boldsymbol{a}_l \Delta^{lk} \middle| \boldsymbol{a}_m \Delta^{mj} \right\rangle}{A^2} =$$

$$= \alpha^{(k)}\alpha^{(j)} \frac{\Delta^{lk}\Delta^{mj}A_{lm}}{A^2} = \alpha^{(k)}\alpha^{(j)} \frac{\Delta^{lk}\delta_l^j A}{A^2} = \alpha^{(k)}\alpha^{(j)} \frac{\Delta^{jk}}{A},$$

$$\boldsymbol{n}^{(k)\,2} = \alpha^{(k)\,2} \frac{\Delta^{kk}}{A}, \quad (j,k = \overline{1, n+m}) \tag{37}$$

Therefore, taking into account formula (18), for non-degenerate hyper-planes we obtain

$$\frac{\left\langle \boldsymbol{n}^{(k)} \middle| \boldsymbol{n}^{(j)} \right\rangle^2}{\boldsymbol{n}^{(k)2}\boldsymbol{n}^{(j)2}} = \frac{(\Delta^{kj})^2}{\Delta^{kk}\Delta^{jj}} = \frac{\left\langle \boldsymbol{a}_1 \wedge \boldsymbol{a}_2 \wedge \ldots \not{\boldsymbol{a}}_k \ldots \wedge \boldsymbol{a}_{n+m} \middle| \boldsymbol{a}_1 \wedge \boldsymbol{a}_2 \wedge \ldots \not{\boldsymbol{a}}_j \ldots \wedge \boldsymbol{a}_{n+m} \right\rangle^2}{\left(\boldsymbol{a}_1 \wedge \boldsymbol{a}_2 \wedge \ldots \not{\boldsymbol{a}}_k \ldots \wedge \boldsymbol{a}_{n+m}\right)^2 \left(\boldsymbol{a}_1 \wedge \boldsymbol{a}_2 \wedge \ldots \not{\boldsymbol{a}}_j \ldots \wedge \boldsymbol{a}_{n+m}\right)^2} \tag{38}$$

(in these formulas summation over repeating indices $j$ and $k$ is not implied).

In 3-dimensional space $\mathrm{span}[\boldsymbol{a}_1, \boldsymbol{a}_2, \boldsymbol{a}_3]$ ($\boldsymbol{a}_1 \wedge \boldsymbol{a}_2 \wedge \boldsymbol{a}_3 \neq 0$) that one has to consider according to Lemma 1, the results obtained give

$$\frac{\left\langle \boldsymbol{n}^{(3)} \middle| \boldsymbol{n}^{(1)} \right\rangle^2}{\boldsymbol{n}^{(3)2}\boldsymbol{n}^{(1)2}} = \frac{(\Delta^{31})^2}{\Delta^{33}\Delta^{11}} = \frac{\left\langle \boldsymbol{a}_1 \wedge \boldsymbol{a}_2 \middle| \boldsymbol{a}_2 \wedge \boldsymbol{a}_3 \right\rangle^2}{(\boldsymbol{a}_1 \wedge \boldsymbol{a}_2)^2 (\boldsymbol{a}_2 \wedge \boldsymbol{a}_3)^2}, \tag{38'}$$

where $\boldsymbol{n}^{(3)}$ and $\boldsymbol{n}^{(1)}$ are normals to non-degenerate 2-planes $\mathrm{span}[\boldsymbol{a}_1, \boldsymbol{a}_2]$ and $\mathrm{span}[\boldsymbol{a}_2, \boldsymbol{a}_3]$ correspondingly.

**Note.** Non-isotropic normal to a degenerate 2-plane is not defined uniquely and belongs to a degenerate 3-dimensional space: if a vector $\boldsymbol{n}$, $\boldsymbol{n}^2 \neq 0$, is orthogonal to a degenerate 2-plane (that contains isotropic vector $\boldsymbol{c}$ and vector $\boldsymbol{a}$, $\boldsymbol{a}^2 \neq 0$, $\langle \boldsymbol{a} | \boldsymbol{c} \rangle = 0$) then, evidently, vectors $\boldsymbol{n} + \alpha\boldsymbol{c}$, $\forall\alpha \in \mathbb{R}$, are orthogonal to the 2-plane as well and, therefore, $\mathrm{span}[\boldsymbol{c}, \boldsymbol{a}, \boldsymbol{n}] \subset \mathbb{P}_{n,m}^{n+m}$, $\boldsymbol{c} \wedge \boldsymbol{a} \wedge \boldsymbol{n} \neq \boldsymbol{0}_{(3)}$, is 3-dimensional degenerate subspace in $\mathbb{P}_{n,m}^{n+m}$ [1]. Nevertheless, the norm of the vector $\boldsymbol{n}$ is defined uniquely:

$$\boldsymbol{c}^2 = 0, \ \langle \boldsymbol{n} | \boldsymbol{c} \rangle = 0 \ \Rightarrow \ (\boldsymbol{n} + \alpha\boldsymbol{c})^2 = \boldsymbol{n}^2, \ \forall\alpha \in \mathbb{R}.$$

Moreover, if the isotropic vector $\boldsymbol{c}$ is common vector of intersecting 2-planes (either degenerate or not), then the LHS of formula (38') is defined uniquely as well.

**Lemma 3.** *For any vectors* $\boldsymbol{a}, \boldsymbol{b} \in \mathbb{P}_{n,m}^{n+m}$, $\boldsymbol{b}^2 \neq 0$, *there exists a vector*

$$\boldsymbol{a}_\perp^{(b)} = \boldsymbol{a} - \boldsymbol{b}\langle \boldsymbol{a} | \boldsymbol{b} \rangle / \boldsymbol{b}^2,$$

*such that*

$$\left\langle \boldsymbol{a}_\perp^{(b)} \middle| \boldsymbol{b} \right\rangle = 0, \quad \boldsymbol{a} \wedge \boldsymbol{b} = \boldsymbol{a}_\perp^{(b)} \wedge \boldsymbol{b}.$$

**Proof.** It is enough to take into account linearity of wedge $\boldsymbol{a} \wedge \boldsymbol{b}$ and obvious equality $k\boldsymbol{b} \wedge \boldsymbol{b} = \boldsymbol{0}_{(2)}$ for $\forall \boldsymbol{b} \in \mathbb{V}(\mathbb{C})$, $k \in \mathbb{C}$. □

**Note.** Lemma 3 applies to bivectors (viewed as vectors in a pseudo-Euclidean space $\Lambda^2 \mathbb{P}_{n,m}^{n+m}$), however, in general case, it does not hold for 2-blades, since they generally do not form a linear space.

**Statement 3.** *If vectors* $\boldsymbol{a}_1, \boldsymbol{a}_2, \boldsymbol{a}_3, \boldsymbol{a}_4 \in \mathbb{P}_{n,m}^{n+m}$ *satisfy the conditions*

$$\boldsymbol{a}_1 \wedge \boldsymbol{a}_2 \wedge \boldsymbol{a}_3 \wedge \boldsymbol{a}_4 = \boldsymbol{0}_{(4)}, \ \boldsymbol{a}_1 \wedge \boldsymbol{a}_2 \neq k\boldsymbol{a}_3 \wedge \boldsymbol{a}_4 \ (k \in \mathbb{R}), \ (\boldsymbol{a}_1 \wedge \boldsymbol{a}_2)^2 \neq 0, \ (\boldsymbol{a}_1 \wedge \boldsymbol{a}_2)^2 \neq 0,$$

*then formulas (27)-(29₁₋₂) for 2-blades* $\boldsymbol{A} = \boldsymbol{a}_1 \wedge \boldsymbol{a}_2$ *and* $\boldsymbol{B} = \boldsymbol{a}_3 \wedge \boldsymbol{a}_4$ *determine corresponding dihedral angles between non-degenerate 2-planes* $\mathrm{span}[\boldsymbol{a}_1, \boldsymbol{a}_2] \subset \mathbb{P}_{n,m}^{n+m}$ *and* $\mathrm{span}[\boldsymbol{a}_3, \boldsymbol{a}_4] \subset \mathbb{P}_{n,m}^{n+m}$.

**Proof.** According formula (18'), the conditions

$$A^2=(a_1\wedge a_2)^2\neq 0,\quad B^2=(a_3\wedge a_4)^2\neq 0$$

indicates that 2-blades $A=a_1\wedge a_2$ and $B=a_3\wedge a_4$ are not isotropic. According to Lemma 1, intersection of correspondent 2-plains – the vector $x$ defined in (32) and (33) – does exist. Let us assume that nonzero coefficients in (32) are $\tilde{\alpha}$ and $\tilde{\lambda}$. If the vector $x$ is not isotropic, then, according to Lemma 3, we find

$$x\wedge a_2=\tilde{\alpha}(a_1\wedge a_2),\quad x\wedge a_4=-\tilde{\lambda}(a_3\wedge a_4)\ \Rightarrow$$

$$A=a_1\wedge a_2=\tilde{\alpha}^{-1}(x\wedge a_2)=\tilde{\alpha}^{-1}(x\wedge a^{(x)}_{2\perp}),\quad B=a_3\wedge a_4=-\tilde{\lambda}^{-1}(x\wedge a_4)=-\tilde{\lambda}^{-1}(x\wedge a^{(x)}_{4\perp}).$$

Inserting $A$ and $B$ in (27)-(29$_{1\text{-}2}$) and calculating inner products of 2-blades according to (18'), one gets

$$\frac{\langle a_1\wedge a_2|a_3\wedge a_4\rangle^2}{(a_1\wedge a_2)^2(a_3\wedge a_4)^2}=\frac{\langle x\wedge a_2|x\wedge a_4\rangle^2}{(x\wedge a_2)^2(x\wedge a_4)^2}=\frac{\left(x^2\langle a^{(x)}_{2\perp}\,|\,a^{(x)}_{4\perp}\rangle\right)^2}{(x^2a^{(x)\,2}_{2\perp})(x^2a^{(x)\,2}_{4\perp})}=\frac{\langle a^{(x)}_{2\perp}\,|\,a^{(x)}_{4\perp}\rangle^2}{a^{(x)\,2}_{2\perp}a^{(x)\,2}_{4\perp}}. \tag{39}$$

Thus, Lemma 2 and formula (38'), as well as the expression (39) generalize – to the case of space $\mathbb{P}^{n+m}_{n,m}$ and subject to the conditions of the Statement – well-known formulas defining the dihedral angle between 2-planes in Euclidean space (see, e.g., [11]). □

**Statement 4.** *If two pseudo-Euclidean 2-planes intersect along isotropic vector* **c***, then the normals to these 2-planes belong to degenerate 2-plane.*

**Proof.** Let $A=c\wedge a,\ B=c\wedge b$ be 2-blades those determine intersecting pseudo-Euclidean 2-planes. Here, according to Lemma 1, without loss of generality

$$c=t+s,\quad a=\alpha_1 t+\beta_1 s,\quad b=\alpha_2 t+\beta_2 s+d, \tag{40}$$

$$t^2=-s^2>0,\ \langle t|s\rangle=\langle d|s\rangle=\langle d|t\rangle=0,\ c^2=0,\ a^2=(\alpha_1^2-\beta_1^2)t^2,\ b^2=(\alpha_2^2-\beta_2^2)t^2+d^2. \tag{41}$$

Evidently,

$$\langle d|c\rangle=\langle d|a\rangle=0,\quad \langle c|a\rangle=(\alpha_1-\beta_1)t^2,\quad \langle c|b\rangle=(\alpha_2-\beta_2)t^2,\quad \langle a|b\rangle=(\alpha_1\alpha_2-\beta_1\beta_2)t^2, \tag{42}$$

$$c\wedge a=(\beta_1-\alpha_1)t\wedge s,\quad c\wedge b=(\beta_2-\alpha_2)t\wedge s+c\wedge d,\quad c\wedge a\wedge b=(\beta_1-\alpha_1)t\wedge s\wedge d\neq \mathbf{0}_{(3)}. \tag{43}$$

Therefore, according to Lemma 2, the normal to span$[c,b]$ is the vector

$$n^{(b,c)}=\alpha\det\begin{bmatrix} c^2 & \langle c|a\rangle & \langle c|b\rangle\\ c & a & b\\ \langle b|c\rangle & \langle b|a\rangle & b^2\end{bmatrix}\det{}^{-1}\begin{bmatrix} c^2 & \langle c|a\rangle & \langle c|b\rangle\\ \langle a|c\rangle & a^2 & \langle a|b\rangle\\ \langle b|c\rangle & \langle b|a\rangle & b^2\end{bmatrix}.$$

Inserting here (40)-(42) and calculating the determinants one gets

$$n^{(b,c)}=\alpha\,\frac{(\alpha_1-\beta_1)t^2}{(\alpha_1-\beta_1)^2t^2s^2d^2}\left(-cd^2+d(\alpha_2-\beta_2)t^2\right)=\frac{\alpha}{\alpha_1-\beta_1}\left(\frac{c}{t^2}-\frac{d}{d^2}(\alpha_2-\beta_2)\right), \tag{44}$$

$$\Rightarrow\quad n^{(b,c)\,2}=\frac{\alpha^2(\alpha_2-\beta_2)^2}{(\alpha_1-\beta_1)^2d^2}, \tag{45}$$

as far $c\wedge a\wedge b\neq \mathbf{0}_{(3)}\Leftrightarrow \alpha_1-\beta_1\neq 0$.

Besides, since the normal $n^{(a,c)}$ to span$[c,a]$ is, evidently, (unique up to norm) vector $d$, expression (44) shows that the following relations are hold:

$$n^{(b,c)}=\lambda n^{(a,c)}+\mu c,\quad \langle n^{(a,c)}|c\rangle=\langle n^{(b,c)}|c\rangle=0,\quad n^{(b,c)}\wedge n^{(a,c)}=\nu c\wedge d\neq \mathbf{0}_{(2)},\ (\lambda,\mu,\nu\in\mathbb{R})$$

that is the vectors $n^{(a,c)}$ and $n^{(b,c)}$ belong to some degenerate 2-plane. □

**Corollary 6.** Calculation of a dihedral angle between pseudo-Euclidean 2-planes intersecting along an isotropic vector, according formula (38'), gives

$$\frac{\left\langle \boldsymbol{n}^{(\boldsymbol{a},\boldsymbol{c})} \middle| \boldsymbol{n}^{(\boldsymbol{b},\boldsymbol{c})} \right\rangle^2}{\boldsymbol{n}^{(\boldsymbol{a},\boldsymbol{c})2}\boldsymbol{n}^{(\boldsymbol{b},\boldsymbol{c})2}} = \frac{\langle \boldsymbol{a}\wedge\boldsymbol{c} | \boldsymbol{b}\wedge\boldsymbol{c}\rangle^2}{(\boldsymbol{a}\wedge\boldsymbol{c})^2(\boldsymbol{b}\wedge\boldsymbol{c})^2} = 1, \tag{46}$$

since, according to (15) in degenerate 2-plane one has

$$\boldsymbol{n}^{(\boldsymbol{a},\boldsymbol{c})2}\boldsymbol{n}^{(\boldsymbol{b},\boldsymbol{c})2} - \left\langle \boldsymbol{n}^{(\boldsymbol{a},\boldsymbol{c})} \middle| \boldsymbol{n}^{(\boldsymbol{b},\boldsymbol{c})} \right\rangle^2 = 0.$$

Calculating independently the RHS of formula (46) one can check it for such a case: since in pseudo-Euclidean 2-planes $\langle \boldsymbol{c}|\boldsymbol{a}\rangle \neq 0$, $\langle \boldsymbol{c}|\boldsymbol{b}\rangle \neq 0$, then one has

$$\frac{\langle \boldsymbol{a}\wedge\boldsymbol{c} | \boldsymbol{b}\wedge\boldsymbol{c}\rangle^2}{(\boldsymbol{a}\wedge\boldsymbol{c})^2(\boldsymbol{b}\wedge\boldsymbol{c})^2} = \frac{\left(\langle \boldsymbol{a}|\boldsymbol{b}\rangle \boldsymbol{c}^2 - \langle \boldsymbol{a}|\boldsymbol{c}\rangle\langle \boldsymbol{c}|\boldsymbol{b}\rangle\right)^2}{\left(\boldsymbol{a}^2\boldsymbol{c}^2 - \langle \boldsymbol{a}|\boldsymbol{c}\rangle^2\right)\left(\boldsymbol{b}^2\boldsymbol{c}^2 - \langle \boldsymbol{b}|\boldsymbol{c}\rangle^2\right)} = \frac{\left(\langle \boldsymbol{a}|\boldsymbol{c}\rangle\langle \boldsymbol{c}|\boldsymbol{b}\rangle\right)^2}{\langle \boldsymbol{a}|\boldsymbol{c}\rangle^2\langle \boldsymbol{b}|\boldsymbol{c}\rangle^2} = 1.$$

Note, that RHS of formula (39) in this case is meaningless: in pseudo-Euclidean 2-planes span$[\boldsymbol{c},\boldsymbol{a}]$ and span$[\boldsymbol{c},\boldsymbol{b}]$ do not exist vectors those are orthogonal to common isotropic vector $\boldsymbol{c}$ of these 2-planes.

**Note**. Degenerate 2-plane can intersect another 2-plane along a non-isotropic vector

$$\boldsymbol{x} = \text{span}[\boldsymbol{c},\boldsymbol{x}] \cap \text{span}[\boldsymbol{a},\boldsymbol{x}],\ \ \boldsymbol{x}^2 \neq 0,\ \ \langle \boldsymbol{c}|\boldsymbol{x}\rangle = \langle \boldsymbol{a}|\boldsymbol{x}\rangle = 0,\ \ \boldsymbol{c} \neq k\boldsymbol{a},$$

or along the isotropic one

$$\boldsymbol{c} = \text{span}[\boldsymbol{c},\boldsymbol{a}] \cap \text{span}[\boldsymbol{c},\boldsymbol{b}],\ \ \langle \boldsymbol{c}|\boldsymbol{a}\rangle = 0,\ \ \boldsymbol{c}^2 = 0,\ \ \boldsymbol{a} \neq k\boldsymbol{b},\ \ k \in \mathbb{R}.$$

Then, for non-isotropic common vector $\boldsymbol{x}$, $\boldsymbol{x}^2 \neq 0$, one gets

$$\langle \boldsymbol{c}\wedge\boldsymbol{x} | \boldsymbol{a}\wedge\boldsymbol{x}\rangle = \langle \boldsymbol{c}|\boldsymbol{a}\rangle \boldsymbol{x}^2 - \langle \boldsymbol{c}|\boldsymbol{x}\rangle\langle \boldsymbol{a}|\boldsymbol{x}\rangle = \langle \boldsymbol{c}|\boldsymbol{a}\rangle \boldsymbol{x}^2,$$

$$(\boldsymbol{c}\wedge\boldsymbol{x})^2 = 0,\ \ \boldsymbol{x}\wedge\boldsymbol{c} \neq \boldsymbol{0}_{(2)},\ \ \boldsymbol{x}\wedge\boldsymbol{a} \neq \boldsymbol{0}_{(2)}.$$

Therefore, in such a case, dihedral angle cannot be determined if $\langle \boldsymbol{c}|\boldsymbol{a}\rangle \neq 0$. If $\langle \boldsymbol{c}|\boldsymbol{a}\rangle = 0$ then the normals $\boldsymbol{n}^{(\boldsymbol{c},\boldsymbol{x})}$ and $\boldsymbol{n}^{(\boldsymbol{a},\boldsymbol{x})}$ to the intersecting 2-planes span$[\boldsymbol{c},\boldsymbol{x}]$ and span$[\boldsymbol{a},\boldsymbol{x}]$, evidently, belong to the degenerate 2-plane span$[\boldsymbol{c},\boldsymbol{a}]$:

$$\boldsymbol{n}^{(\boldsymbol{c},\boldsymbol{x})} = \alpha_1\boldsymbol{c} + \alpha_2\boldsymbol{a},\ \ \boldsymbol{n}^{(\boldsymbol{a},\boldsymbol{x})} = \beta_1\boldsymbol{c} + \beta_2\boldsymbol{a},\ \ \langle \boldsymbol{c}|\boldsymbol{a}\rangle = 0$$

and LHS of formula (38') yield 1, whereas its RHS and both sides of (39) lead to an indeterminacy 0/0.

For the isotropic common vector, one has

$$\langle \boldsymbol{c}\wedge\boldsymbol{a} | \boldsymbol{c}\wedge\boldsymbol{b}\rangle = \boldsymbol{c}^2\langle \boldsymbol{a}|\boldsymbol{b}\rangle - \langle \boldsymbol{c}|\boldsymbol{b}\rangle\langle \boldsymbol{a}|\boldsymbol{c}\rangle = 0,$$

$$(\boldsymbol{c}\wedge\boldsymbol{a})^2 = (\boldsymbol{c}\wedge\boldsymbol{b})^2 = 0,$$

and formula (38') as well as the LHS of formula (39) yield an indeterminacy 0/0 again.

These indeterminacies can be calculated as the corresponding limits (see Appendix 2). Note that although, the vectors orthogonal to the common isotropic vector (along which these 2-planes intersect) are not uniquely determined:

$$\langle \boldsymbol{c}|\boldsymbol{a}\rangle = \langle \boldsymbol{c}|\boldsymbol{b}\rangle = 0 \Rightarrow \langle \boldsymbol{c}|\boldsymbol{a}'\rangle = \langle \boldsymbol{c}|\boldsymbol{b}'\rangle = 0,$$

$$\boldsymbol{a}' = \boldsymbol{a} + \lambda\boldsymbol{c},\ \ \boldsymbol{b}' = \boldsymbol{b} + \mu\boldsymbol{c},\ \ \forall\lambda,\mu \in \mathbb{R},$$

nevertheless

$$\langle \boldsymbol{a}'|\boldsymbol{b}'\rangle = \langle \boldsymbol{a}|\boldsymbol{b}\rangle,\ \ \boldsymbol{a}'^2 = \boldsymbol{a}^2,\ \ \boldsymbol{b}'^2 = \boldsymbol{b}^2.$$

**Corollary 7.** Let us define a vector

$$\boldsymbol{d}=\alpha\boldsymbol{a}+\beta\boldsymbol{b}\in\operatorname{span}[\boldsymbol{a},\boldsymbol{b}]\subset\mathbb{P}_{n,m}^{n+m}$$

such that $\langle\boldsymbol{u}|\boldsymbol{d}\rangle=0$ for some $\boldsymbol{u}\in\mathbb{P}_{n,m}^{n+m}$, $\langle\boldsymbol{u}|\boldsymbol{a}\rangle\neq0$, $\langle\boldsymbol{u}|\boldsymbol{b}\rangle\neq0$. Such vector $\boldsymbol{d}$, evidently, is unique up to (complex) factor:

$$\boldsymbol{d}=k\big(\boldsymbol{a}\langle\boldsymbol{u}|\boldsymbol{b}\rangle-\boldsymbol{b}\langle\boldsymbol{u}|\boldsymbol{a}\rangle\big),\quad k\in\mathbb{C}. \tag{47}$$

Then, using the formula (39), one can define an angle between 2-plane $\operatorname{span}[\boldsymbol{a},\boldsymbol{b}]$ and the vector $\boldsymbol{u}\in\mathbb{P}_{n,m}^{n+m}$, as the dihedral angle between 2-planes $\operatorname{span}[\boldsymbol{a},\boldsymbol{b}]$ and $\operatorname{span}[\boldsymbol{d},\boldsymbol{u}]$. Taking into account (47), (18) and (18'), we obtain (in general complex case; comparison with (39) and other details see in Appendix 3)

$$\langle\boldsymbol{d}\wedge\boldsymbol{u}|\boldsymbol{a}\wedge\boldsymbol{b}\rangle=\langle\boldsymbol{d}|\boldsymbol{a}\rangle\langle\boldsymbol{u}|\boldsymbol{b}\rangle-\langle\boldsymbol{d}|\boldsymbol{b}\rangle\langle\boldsymbol{u}|\boldsymbol{a}\rangle=\big\langle\boldsymbol{d}\big|\boldsymbol{a}\langle\boldsymbol{u}|\boldsymbol{b}\rangle-\boldsymbol{b}\langle\boldsymbol{u}|\boldsymbol{a}\rangle\big\rangle=\boldsymbol{d}^2/k,\quad(\boldsymbol{d}\wedge\boldsymbol{u})^2=\boldsymbol{d}^2\boldsymbol{u}^2, \tag{48}$$

$$\Rightarrow\quad\frac{|\langle\boldsymbol{d}\wedge\boldsymbol{u}|\boldsymbol{a}\wedge\boldsymbol{b}\rangle|^2}{(\boldsymbol{d}\wedge\boldsymbol{u})^2(\boldsymbol{a}\wedge\boldsymbol{b})^2}=\frac{\boldsymbol{d}^2}{|k|^2\boldsymbol{u}^2(\boldsymbol{a}\wedge\boldsymbol{b})^2}=\frac{\big(\boldsymbol{a}\langle\boldsymbol{u}|\boldsymbol{b}\rangle-\boldsymbol{b}\langle\boldsymbol{u}|\boldsymbol{a}\rangle\big)^2}{\boldsymbol{u}^2\big(\boldsymbol{a}^2\boldsymbol{b}^2-|\langle\boldsymbol{a}|\boldsymbol{b}\rangle|^2\big)}=1-\frac{(\boldsymbol{u}\wedge\boldsymbol{a}\wedge\boldsymbol{b})^2}{\boldsymbol{u}^2(\boldsymbol{a}\wedge\boldsymbol{b})^2}. \tag{49}$$

**Note.** The dihedral angles defined above are planar angles in a 2-plane orthogonal to the common vector of the intersecting 2-blades. The type of this 2-plane is entirely determined by the sign of $D$ in expression (10), involving vectors $\boldsymbol{a}$, $\boldsymbol{b}$ belonging to this 2-plane. Correspondingly, for the defined dihedral angle one has to calculate cosh, cos, defect of orthogonality or deviation of collinearity (see (24)-(26)), depending on the sign of $D$.

**5. Some further generalizations and discussion.** In contrast with an Euclidean space, in a pseudo-Euclidean space relations like to CBS inequality imply all three possible cases, those in a case of complex linear space have the forms:

$$|\langle\boldsymbol{a}|\boldsymbol{b}\rangle|^2>\boldsymbol{a}^2\boldsymbol{b}^2,\ \ |\langle\boldsymbol{a}|\boldsymbol{b}\rangle|^2<\boldsymbol{a}^2\boldsymbol{b}^2,\ \ |\langle\boldsymbol{a}|\boldsymbol{b}\rangle|^2=\boldsymbol{a}^2\boldsymbol{b}^2,$$

depending on the kind of the 2-plane (in the complex linear space $\mathbb{P}_{n,m}^{n+m}$) the vectors $\boldsymbol{a},\boldsymbol{b}$ belong to. Correspondingly, an angle and its generalizations between vectors that belong to (real) 2-planes of different kinds are defined, as well as generalizations of dihedral angles between 2-planes (2-blades) of different kinds in real $\mathbb{P}_{n,m}^{n+m}$.

Some authors (see, e. g., [13]) consider a generalization of CBS inequality for number of vectors in a Euclidean space using quadratic forms and Gram's matrix formalism. Evidently, these methods can be applied to pseudo-Euclidean spaces as well (for general case of a complex space, as far the Gram's matrix is Hermitian). Namely, expression (6) that is the source of the CBS inequality (in a complex case) and its generalizations for pseudo-Euclidean space, can be written as quadratic 2-form. In this sense, the CBS inequality in Euclidean space reduces to the conditions for positive definiteness of quadratic 2-forms (see, e. g., [12]). Evidently, our results for pseudo-Euclidean space can be presented as conditions for negative definiteness, positive definiteness and sign-indefiniteness of 2-forms for $\mathbb{P}_{2,0}^2$, $\mathbb{P}_{0,2}^2$ and $\mathbb{P}_{1,1}^2$ spaces, correspondingly. For degenerate spaces one has, obviously (for linearly independent vectors $\boldsymbol{a},\boldsymbol{b}$), non-negative (for ***c-t***-kind space $\operatorname{span}[\boldsymbol{a},\boldsymbol{b}]$), and non-positive 2-forms (for ***c-s***-kind space $\operatorname{span}[\boldsymbol{a},\boldsymbol{b}]$).

From this perspective, the scope of consideration can be naturally extended including in analysis the set of vectors $\{\boldsymbol{a}_k\,|\,\boldsymbol{a}_k\in\mathbb{P}_{n,m}^{n+m}\}$, $1\le k\le n+m$, and assuming that

$$\boldsymbol{u}=\sum\nolimits_{j=1}^{k}\alpha^j\boldsymbol{a}_j\quad(\alpha^j\in\mathbb{C},\ j=\overline{1,k}).$$

Then the expression (6) is transformed and takes the form

$$\langle \boldsymbol{u}|\boldsymbol{u}\rangle = (\bar{\alpha}^1, \bar{\alpha}^2, ..., \bar{\alpha}^k)\left[\langle \boldsymbol{a}_i|\boldsymbol{a}_j\rangle\right]_{i,j=1,k}(\alpha^1, \alpha^2, ..., \alpha^k)^T, \quad (50)$$

where scalar products $\langle \boldsymbol{a}_i|\boldsymbol{a}_j\rangle$, $i, j = \overline{1,k}$, have to be calculated according to (4). By applying Jacobi's theorem regarding the definiteness (positive, negative or indefinite) of quadratic $k$-forms (see, e. g., [12]) one can determine whether the given set of vectors $\{\boldsymbol{a}_i \,|\, i = \overline{1,k}\}$ spans Euclidean subspace $\mathbb{P}^k_{k,0} \subset \mathbb{P}^{m+n}_{m,n}$ (or $\mathbb{P}^k_{0,k} \subset \mathbb{P}^{m+n}_{m,n}$) or pseudo-Euclidean one. Moreover, the number of sign changes in the sequence of leading principal minors of the Gram's matrix (all these minors are real as far the Gramm's matrix is Hermitian) uniquely determines the signature of corresponding subspace in the pseudo-Euclidean space $\mathbb{P}^{m+n}_{m,n}$.

In this context, it is worth noting the following point: the method based on Jacobi's theorem, (on the Gram's matrix method) allows one to draw conclusions about the nature of the subspace to which a given set $\{\boldsymbol{a}_k | \boldsymbol{a}_k \in \mathbb{P}^{n+m}_{n,m}\}$ of $k$ vectors belongs $(1 \le k \le n+m)$, as well as about the nature of the subspaces spanned by its subsets. Moreover, in the general case, one can determine which linear combinations of vectors from a given set yield exclusively time-like vectors, exclusively space-like vectors, vectors of both types, or isotropic vectors. Our results (7)-(12) indicate the coefficients that should be chosen to obtain time-like, space-like, or isotropic vectors (if such exist). It would be of interest to investigate similar problems for sets containing more than two (linearly independent) vectors.

For special case for set of 4 vectors we have obtained similar criterion using CBS-type inequalities for bivectors (see formulas ($23_{1\text{-}3}$)). Thus, for analyzing a span$[\boldsymbol{a}_1, \boldsymbol{a}_2, \boldsymbol{a}_3, \boldsymbol{a}_4]$ we have two different sets of conditions – formulas ($23_{1\text{-}3}$), and Jakobi's conditions for sign-definiteness of the quadratic 4-form (see expression (50) for $k=4$).

For example, if the Gram's matrix of linearly independent vectors $\boldsymbol{a}, \boldsymbol{b}, \boldsymbol{c}, \boldsymbol{d}$ those belong to a (complex) pseudo-Euclidean space has principal leading (real) minors

$$\boldsymbol{a}^2 > 0,\ \det\begin{bmatrix} \boldsymbol{a}^2 & \langle \boldsymbol{a}|\boldsymbol{b}\rangle \\ \langle \boldsymbol{b}|\boldsymbol{a}\rangle & \boldsymbol{b}^2 \end{bmatrix} < 0,\ \det\begin{bmatrix} \boldsymbol{a}^2 & \langle \boldsymbol{a}|\boldsymbol{b}\rangle & \langle \boldsymbol{a}|\boldsymbol{c}\rangle \\ \langle \boldsymbol{b}|\boldsymbol{a}\rangle & \boldsymbol{b}^2 & \langle \boldsymbol{b}|\boldsymbol{c}\rangle \\ \langle \boldsymbol{c}|\boldsymbol{a}\rangle & \langle \boldsymbol{c}|\boldsymbol{b}\rangle & \boldsymbol{c}^2 \end{bmatrix} > 0,\ \det\begin{bmatrix} \boldsymbol{a}^2 & \langle \boldsymbol{a}|\boldsymbol{b}\rangle & \langle \boldsymbol{a}|\boldsymbol{c}\rangle & \langle \boldsymbol{a}|\boldsymbol{d}\rangle \\ \langle \boldsymbol{b}|\boldsymbol{a}\rangle & \boldsymbol{b}^2 & \langle \boldsymbol{b}|\boldsymbol{c}\rangle & \langle \boldsymbol{b}|\boldsymbol{d}\rangle \\ \langle \boldsymbol{c}|\boldsymbol{a}\rangle & \langle \boldsymbol{c}|\boldsymbol{b}\rangle & \boldsymbol{c}^2 & \langle \boldsymbol{c}|\boldsymbol{d}\rangle \\ \langle \boldsymbol{d}|\boldsymbol{a}\rangle & \langle \boldsymbol{d}|\boldsymbol{b}\rangle & \langle \boldsymbol{d}|\boldsymbol{c}\rangle & \boldsymbol{d}^2 \end{bmatrix} < 0,$$

then the number of sign changes in this sequence is 3, and one has to conclude that $\boldsymbol{a}, \boldsymbol{b}, \boldsymbol{c}, \boldsymbol{d} \in \mathbb{P}^4_{1,3}$. Note, that formula ($23_1$) holds for the case when pseudo-Euclidean space is span$[\boldsymbol{a} \wedge \boldsymbol{b}, \boldsymbol{c} \wedge \boldsymbol{d}]$ with space-like 2-blades $\boldsymbol{a} \wedge \boldsymbol{b}$ and $\boldsymbol{c} \wedge \boldsymbol{d}$ $((\boldsymbol{a} \wedge \boldsymbol{b})^2 < 0, (\boldsymbol{c} \wedge \boldsymbol{d})^2 < 0)$. Apparently, the Jacobi method allows the space to be described in greater detail. It seems interesting to investigate what kind of restrictions give such sets of inequalities for $k > 4$.

Another interesting problem is to investigate possibility of generalization for a pseudo-Euclidean space inequalities derived in [14], [15] for 4 vectors in Euclidean space. From this perspective, it should be noted that [15] considers sequences that can be viewed as sets of components of vectors – elements of a linear space – only in the real case; in the complex case, however, such sequences do not constitute vector components. Consequently, one cannot speak of generalizing the results obtained in [15] to an arbitrary (complex) linear space – in particular, to a pseudo-Euclidean space.

For the "angle" between vectors $\boldsymbol{a}, \boldsymbol{b}$ belonging to a complex Euclidean space, some authors propose the formula

$$\cos\varphi_{(\boldsymbol{a},\boldsymbol{b})} = \frac{\mathrm{Re}\langle \boldsymbol{a}|\boldsymbol{b}\rangle}{|\boldsymbol{a}||\boldsymbol{b}|} = \frac{\left(\langle \boldsymbol{a}|\boldsymbol{b}\rangle + \langle \boldsymbol{b}|\boldsymbol{a}\rangle\right)}{2|\boldsymbol{a}||\boldsymbol{b}|}, \tag{51}$$

assuming that for these vectors the following inequality holds

$$\left|\mathrm{Re}\langle \boldsymbol{a}|\boldsymbol{b}\rangle\right|^2 - \boldsymbol{a}^2\boldsymbol{b}^2 \le 0. \tag{52}$$

For real case formulas (51) and (52) evidently reduce to (24) and to CBS inequality, correspondingly. However, we consider it important to note some incorrectness of these expressions. Namely, in the case $\boldsymbol{b} = \lambda\boldsymbol{a},\ \lambda \in \mathbb{C},$ Theorem 2 proved above gives $\boldsymbol{a}^2\boldsymbol{b}^2 = \left|\langle \boldsymbol{a}|\boldsymbol{b}\rangle\right|^2$ (which must also hold in the special case of Euclidean space) while the LHS of expression (52) in this case gives

$$\left|\mathrm{Re}\langle \boldsymbol{a}|\lambda\boldsymbol{a}\rangle\right|^2 - \boldsymbol{a}^2\langle \lambda\boldsymbol{a}|\lambda\boldsymbol{a}\rangle = |\mathrm{Re}\lambda|^2\langle \boldsymbol{a}|\boldsymbol{a}\rangle^2 - |\lambda|^2\boldsymbol{a}^2\langle \boldsymbol{a}|\boldsymbol{a}\rangle = \langle \boldsymbol{a}|\boldsymbol{a}\rangle^2\left(|\mathrm{Re}\lambda|^2 - |\lambda|^2\right),$$

that, in general, is negative (being 0 for real $\lambda$ only). Hence, expression

$$\left|\mathrm{Re}\langle \boldsymbol{a}|\boldsymbol{b}\rangle\right|^2 - \boldsymbol{a}^2\boldsymbol{b}^2 \lesseqgtr 0$$

can be considered as alternative form of CBS inequality and its generalization only for real spaces. Correspondingly, one cannot define angles between vectors in complex spaces based on formula (52).

Besides, in our view, the definitions of the *Lorentzian time-like angle* $\delta_{(\boldsymbol{a},\boldsymbol{b})}$ introduced by some authors (see, e.g., [4], [5]) is wrong, due to for that definition one has $\sinh\delta_{(\boldsymbol{a},\boldsymbol{b})} = 0 \Leftrightarrow \delta_{(\boldsymbol{a},\boldsymbol{b})} = 0$ iff $\langle \boldsymbol{a}|\boldsymbol{b}\rangle = 0,$ that is, iff vectors $\boldsymbol{a}, \boldsymbol{b} \in \mathbb{P}^{n+m}_{n,m}$ are orthogonal. For this reason, we call $\delta_{(\boldsymbol{a},\boldsymbol{b})}$ defined in ($25_1$) *defect of orthogonality*.

If in a real pseudo-Euclidean space $\mathbb{P}^{n+m}_{n,m}$ one has $n=1$ (or $m=1$) we get "Minkovsky-like" space (in Minkovsky pseudo-Euclidean space that is using in physics one has $n=1\ m=3$). In contrast with a general case, in "Minkovsky-like" spaces exists only one kind of degenerate 2-planes – $\boldsymbol{c}$-$\boldsymbol{s}$- (or $\boldsymbol{c}$-$\boldsymbol{t}$-) kind only; $\boldsymbol{c}$-$\boldsymbol{c}$-kind degenerate 2-planes in such spaces do not exist.

Besides, in "Minkovsky-like" space $\mathbb{P}^{1+m}_{1,m}$ we have $\boldsymbol{t}^2 > 0,\ \langle \boldsymbol{t}|\boldsymbol{d}\rangle = 0 \Rightarrow \boldsymbol{d}^2 < 0.$ Therefore, taking into consideration that (see formula (A3.6) in Appendix 3)

$$\boldsymbol{t}_{\|}^{(\boldsymbol{a},\boldsymbol{b})\,2}(\boldsymbol{a}\wedge\boldsymbol{b})^2 = \left(\boldsymbol{a}\langle \boldsymbol{t}|\boldsymbol{b}\rangle - \boldsymbol{b}\langle \boldsymbol{t}|\boldsymbol{a}\rangle\right)^2, \quad \boldsymbol{a},\boldsymbol{b} \in \mathbb{P}^{1+m}_{1,m},\ k \in \mathbb{C},$$

one gets (due to the vector $\boldsymbol{a}\langle \boldsymbol{t}|\boldsymbol{b}\rangle - \boldsymbol{b}\langle \boldsymbol{t}|\boldsymbol{a}\rangle,$ evidently, is orthogonal to the vector $\boldsymbol{t}$)

$$\boldsymbol{t}_{\|}^{(\boldsymbol{a},\boldsymbol{b})\,2}(\boldsymbol{a}\wedge\boldsymbol{b})^2 < 0. \tag{53}$$

It means that in "Minkovsky-like" spaces $\mathbb{P}^{1+m}_{1,m}$ the sign of the expression $(\boldsymbol{a}\wedge\boldsymbol{b})^2 = D$ (that determines the direction in the CBS-type inequality) is opposite to the sign of $\boldsymbol{t}_{\|}^{(\boldsymbol{a},\boldsymbol{b})\,2}$:

$$D<0 \ \Leftrightarrow\ \boldsymbol{t}_{\|}^{(\boldsymbol{a},\boldsymbol{b})\,2} > 0; \quad D>0 \ \Leftrightarrow\ \boldsymbol{t}_{\|}^{(\boldsymbol{a},\boldsymbol{b})\,2} < 0.$$

Given that the span$[\boldsymbol{t}, \boldsymbol{t}_{\|}^{(\boldsymbol{a},\boldsymbol{b})}]$ intersecting the span$[\boldsymbol{a},\boldsymbol{b}]$ along the vector $\boldsymbol{t}_{\|}^{(\boldsymbol{a},\boldsymbol{b})}$ is evidently pseudo-Euclidean, one can conclude that if the "angle" between the axis $\boldsymbol{t}$ and 2-plane span$[\boldsymbol{a},\boldsymbol{b}]$ – deviation of collinearity, see ($25_2$) – is "small enough" (cosh of this angle is finite) then the 2-plane span$[\boldsymbol{a},\boldsymbol{b}]$ is pseudo-Euclidean; and if the 2-plane $\mathrm{span}[\boldsymbol{a},\boldsymbol{b}]$ is Euclidean and, therefore, if $\boldsymbol{t}_{\|}^{(\boldsymbol{a},\boldsymbol{b})\,2} < 0$, then the "angle" between the $\boldsymbol{t}$ axis and the 2-plane $\mathrm{span}[\boldsymbol{a},\boldsymbol{b}]$ characterizes the "defect of orthogonality" $\delta$ (see ($25_1$)):

$$\delta = (-\boldsymbol{t}_{\|}^{(\boldsymbol{a},\boldsymbol{b})\,2}\boldsymbol{t}^2)^{-1/2}\sinh^{-1}\left\langle \boldsymbol{t}_{\|}^{(\boldsymbol{a},\boldsymbol{b})} \middle| \boldsymbol{t} \right\rangle = \sinh^{-1}(-\boldsymbol{t}_{\|}^{(\boldsymbol{a},\boldsymbol{b})\,2})^{1/2}. \quad (\boldsymbol{t}^2 = 1)$$

If $D \to 0_{\pm}$, that is if the 2-plane $\mathrm{span}[\boldsymbol{a},\boldsymbol{b}]$ becomes degenerate ($\boldsymbol{s}$-$\boldsymbol{c}$-kind), then, according to (53), $\boldsymbol{t}_{\|}^{(\boldsymbol{a},\boldsymbol{b})\,2} \to \mp\infty$. In the case $D \to 0_{-}$, $\boldsymbol{t}_{\|}^{(\boldsymbol{a},\boldsymbol{b})\,2} \to +\infty$ we get

$$\cosh\varphi = (\boldsymbol{t}_{\|}^{(\boldsymbol{a},\boldsymbol{b})\,2})^{1/2} \to +\infty,$$

and in the case $D \to 0_{+}$, $\boldsymbol{t}_{\|}^{(\boldsymbol{a},\boldsymbol{b})\,2} \to -\infty$, correspondingly –

$$\sinh\delta = (-\boldsymbol{t}_{\|}^{(\boldsymbol{a},\boldsymbol{b})\,2})^{1/2} \to +\infty \quad (\boldsymbol{t}^2 = 1).$$

It should be noted that our definition ($25_1$) of a "defect of orthogonality" $\delta_{(\boldsymbol{a},\boldsymbol{b})}$ as well as definition ($25_2$) of a "deviation of collinearity" $\varphi_{(\boldsymbol{a},\boldsymbol{b})}$ does not depend on a choice of basic vectors $\boldsymbol{t}$ and $\boldsymbol{s}$ in the 2-plane (in the pseudo-Euclidean space $\mathbb{P}_{n,m}^{n+m}$) the (real) vectors $\boldsymbol{a},\boldsymbol{b} \in \mathbb{P}_{n,m}^{n+m}$ belong to.

## R E F E R E N C E S

**Appendix 1.**

The equation (6) and expressions (8)-(10) give

$$\boldsymbol{u}^2(z) = |\alpha|^2\boldsymbol{b}^2\Big(|\,z + \langle\boldsymbol{b}|\boldsymbol{a}\rangle/\boldsymbol{b}^2\,|^2 - D/\langle\boldsymbol{b}|\boldsymbol{b}\rangle^2\Big) = 0 \qquad \text{(A1.1)}$$

for $z = z_0 + Re^{i\varphi}$, $R > 0$, $\varphi \in (-\pi;\pi]$. So, if $D > 0$, $\boldsymbol{b}^2 \neq 0$, there exists continuous set of isotropic vectors $\boldsymbol{u}(z) \neq \boldsymbol{0}$. But only 2 vectors of this set are linearly independent, not of them being mutually orthogonal.

Indeed, choosing $z_1 = z_0 + Re^{i\varphi_1}$, $z_2 = z_0 + Re^{i\varphi_2}$, one gets linearly independent vectors $\boldsymbol{u}(z_1)$ and $\boldsymbol{u}(z_2)$:

$$\boldsymbol{u}_k = \boldsymbol{u}(z_k) = \alpha(\boldsymbol{a} + z_k\boldsymbol{b}) = \alpha(\boldsymbol{a} + z_0\boldsymbol{b} + Re^{i\varphi_k}\boldsymbol{b}), \quad k = 1,2.$$

$$A\boldsymbol{u}_1 + B\boldsymbol{u}_2 = \boldsymbol{0} \;\Rightarrow\; \alpha\{(A+B)\boldsymbol{a} + (z_1A + z_2B)\boldsymbol{b}\} = \boldsymbol{0} \;\; (\boldsymbol{a}\wedge\boldsymbol{b} \neq \boldsymbol{0}_{(2)},\; A, B \in \mathbb{C}) \;\Rightarrow\; \begin{cases} A + B = 0, \\ z_1A + z_2B = 0, \end{cases}$$

$$\det\begin{bmatrix} 1 & 1 \\ z_1 & z_2 \end{bmatrix} \equiv \Delta = z_2 - z_1 = R(e^{i\varphi_2} - e^{i\varphi_1}) = Re^{i\varphi_1}(e^{i(\varphi_2-\varphi_1)} - 1) = 0 \;\Leftrightarrow\; \varphi_2 = \varphi_1 \quad (\varphi_1, \varphi_2 \in (-\pi;\pi])$$

$$\Rightarrow\; \varphi_2 \neq \varphi_1,\; \Delta \neq 0 \;\Rightarrow\; A = B = 0.$$

Besides, these vectors are not mutually orthogonal:

$$\langle\boldsymbol{u}_1|\boldsymbol{u}_2\rangle = |\alpha|^2\Big(\boldsymbol{a}^2 + \langle\boldsymbol{a}|\boldsymbol{b}\rangle z_2 + \langle\boldsymbol{b}|\boldsymbol{a}\rangle\bar{z}_1 + \boldsymbol{b}^2\bar{z}_1z_2\Big) = |\alpha|^2/\boldsymbol{b}^2\Big(\boldsymbol{a}^2\boldsymbol{b}^2 - |\langle\boldsymbol{a}|\boldsymbol{b}\rangle|^2 + De^{i(\varphi_2-\varphi_1)}\Big) =$$

$$= |\alpha|^2 D/\boldsymbol{b}^2\Big(e^{i(\varphi_2-\varphi_1)} - 1\Big) \neq 0 \;\Leftrightarrow\; \varphi_2 \neq \varphi_1. \;\; (D > 0,\; \boldsymbol{b}^2 \neq 0)$$

**Note**. In contrast with a real case where coefficients of a quadratic equation are completely determined by 2 roots of the equation (according to the Viet's Theorem), to determine the coefficients of the equation (A1.1) it is necessary (and sufficient) to know 3 roots of it.

Indeed, assuming that $z_k = x_k + iy_k$, $k = 1,2,3$, $z_0 = x_0 + iy_0$, one has

$$|z_k - z_0|^2 - R^2 = 0 \;\Leftrightarrow\; (x_k - x_0)^2 + (y_k - y_0)^2 = R^2 \;\Rightarrow \qquad \text{(A1.2)}$$

$$x_0^2 - 2x_0x_k + y_0^2 - 2y_0y_k = R^2 - x_k^2 - y_k^2, \quad k = 1,2,3. \qquad \text{(A1.3)}$$

As far (A1.3) is a system of equations in three real variables $x_0, y_0, R$ (those determine coefficients of the equation (A1.1)) and as far $x_k, y_k$, $k = 1,2,3$, are known, we can reduce it to a system of linear equations

$$\begin{cases} 2x_0(x_2 - x_1) + 2y_0(y_2 - y_1) = |z_2|^2 - |z_1|^2, \\ 2x_0(x_3 - x_1) + 2y_0(y_3 - y_1) = |z_3|^2 - |z_1|^2, \end{cases}$$

that has solutions

$$x_0 = \frac{\det\begin{bmatrix} |z_2|^2 - |z_1|^2 & y_2 - y_1 \\ |z_3|^2 - |z_1|^2 & y_3 - y_1 \end{bmatrix}}{2\det\begin{bmatrix} x_2 - x_1 & y_2 - y_1 \\ x_3 - x_1 & y_3 - y_1 \end{bmatrix}} = \frac{\det\begin{bmatrix} |z_1|^2 & y_1 & 1 \\ |z_2|^2 & y_2 & 1 \\ |z_3|^2 & y_3 & 1 \end{bmatrix}}{2\det\begin{bmatrix} x_1 & y_1 & 1 \\ x_2 & y_2 & 1 \\ x_3 & y_3 & 1 \end{bmatrix}}, \quad y_0 = \frac{\det\begin{bmatrix} x_1 & |z_1|^2 & 1 \\ x_2 & |z_2|^2 & 1 \\ x_3 & |z_3|^2 & 1 \end{bmatrix}}{2\det\begin{bmatrix} x_1 & y_1 & 1 \\ x_2 & y_2 & 1 \\ x_3 & y_3 & 1 \end{bmatrix}},$$

and then find $R$ (that is real positive) from (A1.2).

Evidently, the result obtained is consistent with the well-known fact that three points determine a unique circle (see, e.g., [11], formula (2.5-3)).

**Appendix 2.**

*a*) According to Lemma 1, let us consider 3-dimensional pseudo-Euclidean space

$$\mathrm{span}[\boldsymbol{a},\boldsymbol{b},\boldsymbol{c}]=\{\boldsymbol{u}\,|\,\boldsymbol{u}=\alpha\boldsymbol{a}+\beta\boldsymbol{b}+\gamma\boldsymbol{c},\ \boldsymbol{c}^2=0,\ \langle\boldsymbol{a}|\boldsymbol{c}\rangle=\langle\boldsymbol{b}|\boldsymbol{c}\rangle=0,\ \boldsymbol{a},\boldsymbol{b},\boldsymbol{c}\neq\boldsymbol{0},\ \alpha,\beta,\gamma\in\mathbb{R}\}\subset\mathbb{P}^{n+m}_{n,m}$$

and let $\mathrm{span}[\boldsymbol{c},\boldsymbol{a}]$ and $\mathrm{span}[\boldsymbol{c},\boldsymbol{b}]$ be intersecting degenerate 2-planes in this space (evidently, the vectors $\boldsymbol{a},\boldsymbol{b}\in\mathbb{P}^{n+m}_{n,m}$ cannot be isotropic, $\boldsymbol{a}^2\neq 0,\ \boldsymbol{b}^2\neq 0,$ since there are no $\boldsymbol{c}$-$\boldsymbol{c}$- kind degenerate 2-planes in a 3-dimensional pseudo-Euclidean space, e. (see, e. g., [1]]); moreover, the space $\mathbb{P}^{n+m}_{n,m}$ that is being under consideration has to have dimensions $n\geq 1,\ m\geq 3\ \Rightarrow\ n+m\geq 4$).

Let

$$\varepsilon=\lambda\boldsymbol{a}+\mu\boldsymbol{b}+\nu\boldsymbol{c}\in\mathrm{span}[\boldsymbol{a},\boldsymbol{b},\boldsymbol{c}],\ \lambda,\mu\neq 0\ \Rightarrow$$

$$\langle\varepsilon|\boldsymbol{a}\rangle=\lambda\boldsymbol{a}^2+\mu\langle\boldsymbol{a}|\boldsymbol{b}\rangle,\ \langle\varepsilon|\boldsymbol{b}\rangle=\mu\boldsymbol{b}^2+\lambda\langle\boldsymbol{a}|\boldsymbol{b}\rangle,\ \langle\varepsilon|\boldsymbol{c}\rangle=0;\ \ \varepsilon^2=\lambda^2\boldsymbol{a}^2+\mu^2\boldsymbol{b}^2+2\lambda\mu\langle\boldsymbol{a}|\boldsymbol{b}\rangle \tag{A2.1}$$

and calculate a limit

$$\frac{\langle\boldsymbol{c}\wedge\boldsymbol{a}\,|\,\boldsymbol{c}\wedge\boldsymbol{b}\rangle^2}{(\boldsymbol{c}\wedge\boldsymbol{a})^2(\boldsymbol{c}\wedge\boldsymbol{b})^2}=\lim_{\varepsilon\to 0}\frac{\langle(\boldsymbol{c}+\varepsilon)\wedge\boldsymbol{a}\,|\,(\boldsymbol{c}+\varepsilon)\wedge\boldsymbol{b}\rangle^2}{\big((\boldsymbol{c}+\varepsilon)\wedge\boldsymbol{a}\big)^2\big((\boldsymbol{c}+\varepsilon)\wedge\boldsymbol{b}\big)^2}. \tag{A2.2}$$

According to (18') and (A2.1) we have

$$\begin{aligned}
&\langle(\boldsymbol{c}+\varepsilon)\wedge\boldsymbol{a}\,|\,(\boldsymbol{c}+\varepsilon)\wedge\boldsymbol{b}\rangle=(\boldsymbol{c}+\varepsilon)^2\langle\boldsymbol{a}|\boldsymbol{b}\rangle-\langle\boldsymbol{c}+\varepsilon|\boldsymbol{b}\rangle\langle\boldsymbol{a}|\boldsymbol{c}+\varepsilon\rangle=\varepsilon^2\langle\boldsymbol{a}|\boldsymbol{b}\rangle-\langle\varepsilon|\boldsymbol{b}\rangle\langle\boldsymbol{a}|\boldsymbol{c}+\varepsilon\rangle=\\
&\quad=\big(\lambda^2\boldsymbol{a}^2+\mu^2\boldsymbol{b}^2+2\lambda\mu\langle\boldsymbol{a}|\boldsymbol{b}\rangle\big)\langle\boldsymbol{a}|\boldsymbol{b}\rangle-\big(\lambda\boldsymbol{a}^2+\mu\langle\boldsymbol{a}|\boldsymbol{b}\rangle\big)\big(\mu\boldsymbol{b}^2+\lambda\langle\boldsymbol{a}|\boldsymbol{b}\rangle\big)=-\lambda\mu D,\\
&\big((\boldsymbol{c}+\varepsilon)\wedge\boldsymbol{a}\big)^2=(\boldsymbol{c}+\varepsilon)^2\boldsymbol{a}^2-\langle\boldsymbol{c}+\varepsilon|\boldsymbol{a}\rangle^2=\varepsilon^2\boldsymbol{a}^2-\langle\varepsilon|\boldsymbol{a}\rangle^2=\boldsymbol{a}^2\langle\varepsilon|\varepsilon_\perp^{(\boldsymbol{a})}\rangle=\\
&\quad=\big(\lambda^2\boldsymbol{a}^2+\mu^2\boldsymbol{b}^2+2\lambda\mu\langle\boldsymbol{a}|\boldsymbol{b}\rangle\big)\boldsymbol{a}^2-\big(\lambda\boldsymbol{a}^2+\mu\langle\boldsymbol{a}|\boldsymbol{b}\rangle\big)^2=\mu^2 D,\\
&\big((\boldsymbol{c}+\varepsilon)\wedge\boldsymbol{b}\big)^2=\lambda^2 D,
\end{aligned}$$

where

$$D=\boldsymbol{a}^2\boldsymbol{b}^2-\langle\boldsymbol{a}|\boldsymbol{b}\rangle^2\neq 0, \tag{A2.3}$$

as far the $\mathrm{span}[\boldsymbol{a},\boldsymbol{b}]$ cannot be degenerate 2-plane: in 3-dimensional pseudo-Euclidean space $\mathrm{span}[\boldsymbol{a},\boldsymbol{b},\boldsymbol{c}]$ there exists no isotropic vector that belongs to $\mathrm{span}[\boldsymbol{a},\boldsymbol{b}]$, $\langle\boldsymbol{a}|\boldsymbol{c}\rangle=\langle\boldsymbol{b}|\boldsymbol{c}\rangle=0$, and therefore is orthogonal to the isotropic vector $\boldsymbol{c}$.

Inserting in (A2.2) the results obtained we get

$$\frac{\langle\boldsymbol{c}\wedge\boldsymbol{a}\,|\,\boldsymbol{c}\wedge\boldsymbol{b}\rangle^2}{(\boldsymbol{c}\wedge\boldsymbol{a})^2(\boldsymbol{c}\wedge\boldsymbol{b})^2}=\lim_{\lambda,\mu\to 0}\frac{\lambda^2\mu^2D^2}{\lambda^2D\mu^2D}=1 \tag{A2.4}$$

and, therefore, the dihedral angle between degenerate 2-planes intersecting along common isotropic vector can be defined as correspondent limit in ($29_{1,2}$) (compare with limit of ($26_{1,2}$)):

$$\cos\varphi=\frac{|\langle\boldsymbol{c}\wedge\boldsymbol{a}\,|\,\boldsymbol{c}\wedge\boldsymbol{b}\rangle|}{\sqrt{(\boldsymbol{c}\wedge\boldsymbol{a})^2(\boldsymbol{c}\wedge\boldsymbol{b})^2}}=1, \tag{A2.5}$$

as far both bivectors (2-blades) $\boldsymbol{c}\wedge\boldsymbol{a}$ and $\boldsymbol{c}\wedge\boldsymbol{b}$ are isotropic ones:

$$(\boldsymbol{c}\wedge\boldsymbol{a})^2=\boldsymbol{c}^2\boldsymbol{a}^2-\langle\boldsymbol{c}|\boldsymbol{a}\rangle^2=0,\quad(\boldsymbol{c}\wedge\boldsymbol{b})^2=\boldsymbol{c}^2\boldsymbol{b}^2-\langle\boldsymbol{c}|\boldsymbol{b}\rangle^2=0.$$

However, in the RHS of formula (39) one has $\langle\boldsymbol{a}|\boldsymbol{b}\rangle^2/\boldsymbol{a}^2\boldsymbol{b}^2\neq 1$; consequently, formula (39) itself does not hold in the case under consideration.

*b*) Let now the vector $\boldsymbol{x}\in\mathbb{P}_{n,m}^{n+m}$ common to intersecting 2-planes be not isotropic, $\boldsymbol{x}^2\neq 0$.

*b*1) If one of the intersecting 2-planes is degenerate, say, $\boldsymbol{c}\wedge\boldsymbol{x}, \langle\boldsymbol{c}|\boldsymbol{x}\rangle=0$, while the other, $\boldsymbol{a}\wedge\boldsymbol{x}$, is not (being Euclidean or pseudo-Euclidean), then the numerator in the LHS of (39) becomes $\langle\boldsymbol{c}\,|\,\boldsymbol{a}_{\perp}^{(x)}\rangle^2\boldsymbol{x}^4$, whereas the denominator is zero: $\boldsymbol{c}^2\boldsymbol{x}^4\boldsymbol{a}_{\perp}^{(x)2}=0$. Thus, the only case in which expression (39) can make sense is the case $\langle\boldsymbol{c}|\,\boldsymbol{a}_{\perp}^{(x)}\rangle=0$, that is the case where the 2-plane $\boldsymbol{c}\wedge\boldsymbol{a}_{\perp}^{(x)}$ – the 2-plane in which the dihedral angle is to be measured – is degenerate.

To calculate the limit $\boldsymbol{c}\mapsto\lim\limits_{\boldsymbol{\varepsilon}\to\boldsymbol{0}}(\boldsymbol{c}+\boldsymbol{\varepsilon})$ for such a case (where $\forall\boldsymbol{\varepsilon}\in\mathbb{P}_{n,m}^{n+m}$), preliminary we evaluate the LHS of (39) (bellow, for sake of brevity, we denote $\boldsymbol{a}_{\perp}^{(x)}=\tilde{\boldsymbol{a}}$):

$$\frac{\langle\boldsymbol{x}\wedge(\boldsymbol{c}+\boldsymbol{\varepsilon})\,|\,\boldsymbol{x}\wedge\tilde{\boldsymbol{a}}\rangle^2}{\left(\boldsymbol{x}\wedge(\boldsymbol{c}+\boldsymbol{\varepsilon})\right)^2\left(\boldsymbol{x}\wedge\tilde{\boldsymbol{a}}\right)^2}=\frac{\boldsymbol{x}^2\langle\boldsymbol{c}+\boldsymbol{\varepsilon}\,|\,\tilde{\boldsymbol{a}}\rangle^2}{\left[\boldsymbol{x}^2(\boldsymbol{c}+\boldsymbol{\varepsilon})^2-\langle\boldsymbol{c}+\boldsymbol{\varepsilon}\,|\,\boldsymbol{x}\rangle^2\right]\tilde{\boldsymbol{a}}^2}=\frac{\boldsymbol{x}^2\langle\boldsymbol{\varepsilon}\,|\,\tilde{\boldsymbol{a}}\rangle^2}{\left(\boldsymbol{x}^2\boldsymbol{\varepsilon}^2+2\boldsymbol{x}^2\langle\boldsymbol{c}|\boldsymbol{\varepsilon}\rangle-\langle\boldsymbol{\varepsilon}|\boldsymbol{x}\rangle^2\right)\tilde{\boldsymbol{a}}^2}.$$

According to Lemma 1 any 4 vectors that form 2 intersecting 2-planes belong to a 3-dimensional space. In the case under consideration this space is spanned by (linearly independent) vectors $\boldsymbol{x},\tilde{\boldsymbol{a}},\boldsymbol{c}$. Therefore, let us present $\boldsymbol{\varepsilon}\in\mathbb{P}_{n,m}^{n+m}$ as

$$\boldsymbol{\varepsilon}=\varepsilon_1\tilde{\boldsymbol{a}}+\varepsilon_2\boldsymbol{x}+\varepsilon_3\boldsymbol{c}\ \Rightarrow$$
$$\langle\boldsymbol{\varepsilon}|\tilde{\boldsymbol{a}}\rangle=\varepsilon_1\tilde{\boldsymbol{a}}^2,\quad \langle\boldsymbol{\varepsilon}|\boldsymbol{x}\rangle=\varepsilon_2\boldsymbol{x}^2,\quad \boldsymbol{\varepsilon}^2=\varepsilon_1^2\tilde{\boldsymbol{a}}^2+\varepsilon_2^2\boldsymbol{x}^2,$$
$$\langle\boldsymbol{\varepsilon}|\boldsymbol{c}\rangle=\varepsilon_1\langle\tilde{\boldsymbol{a}}|\boldsymbol{c}\rangle+\varepsilon_2\langle\boldsymbol{x}|\boldsymbol{c}\rangle=0.$$

Here we have taken into account that $\langle\boldsymbol{c}|\tilde{\boldsymbol{a}}\rangle=\langle\boldsymbol{x}|\tilde{\boldsymbol{a}}\rangle=\langle\boldsymbol{c}|\boldsymbol{x}\rangle=0,\ \boldsymbol{c}^2=0$. Then

$$\frac{\boldsymbol{x}^2\langle\boldsymbol{\varepsilon}|\tilde{\boldsymbol{a}}\rangle^2}{\left(\boldsymbol{x}^2\boldsymbol{\varepsilon}^2+2\boldsymbol{x}^2\langle\boldsymbol{c}|\boldsymbol{\varepsilon}\rangle-\langle\boldsymbol{\varepsilon}|\boldsymbol{x}\rangle^2\right)\tilde{\boldsymbol{a}}^2}=\frac{\varepsilon_1^2\tilde{\boldsymbol{a}}^2}{\varepsilon_1^2\tilde{\boldsymbol{a}}^2+\varepsilon_2^2\boldsymbol{x}^2-\varepsilon_2^2\boldsymbol{x}^2}=1.$$

Note, that equality (39) makes sense only if $\langle\boldsymbol{c}+\boldsymbol{\varepsilon}|\boldsymbol{x}\rangle=0\Leftrightarrow\varepsilon_2\boldsymbol{x}^2=0\Leftrightarrow\varepsilon_2=0$, which does not affect the limit. Therefore, for such 2-planes the dihedral angle can be calculated as a correspondent limit of the plane "angle" in the degenerate 2-plane span$[\boldsymbol{c},\tilde{\boldsymbol{a}}]$:

$$\lim_{\boldsymbol{\varepsilon}\to\boldsymbol{0}}\frac{\langle\boldsymbol{c}+\boldsymbol{\varepsilon}|\tilde{\boldsymbol{a}}\rangle^2}{(\boldsymbol{c}+\boldsymbol{\varepsilon})^2\tilde{\boldsymbol{a}}^2}=\lim_{\boldsymbol{\varepsilon}\to\boldsymbol{0}}\frac{\langle\boldsymbol{\varepsilon}|\tilde{\boldsymbol{a}}\rangle^2}{\left(2\langle\boldsymbol{c}|\boldsymbol{\varepsilon}\rangle+\boldsymbol{\varepsilon}^2\right)\tilde{\boldsymbol{a}}^2}=\lim_{\varepsilon_1\to 0}\frac{\varepsilon_1^2\tilde{\boldsymbol{a}}^2}{\varepsilon_1^2\tilde{\boldsymbol{a}}^2}=1.$$

Thus, in this case equality (39) is correct as well.

*b*2) If alternatively, both intersecting 2-planes are degenerate, then the numerator in RHS of (39) is $\langle\boldsymbol{c}_1|\boldsymbol{c}_2\rangle^2$, whereas the denominator is $\boldsymbol{c}_1^2\boldsymbol{c}_2^2=0$. Thus, (39) would retain its meaning only in limit $\langle\boldsymbol{c}_1|\boldsymbol{c}_2\rangle\to 0$. Such limit exists only in spaces $\mathbb{P}_{n,m}^{n+m}$ with dimensions $n,m\geq 2\Rightarrow n+m\geq 4$. Therefore, according to Lemma 1, there does not exist degenerate 2-planes intersecting along non-isotropic vector with finite "dihedral angle" (cosh of correspondent plane angle is infinite).

Note, that according to Lemma 1, none of the intersecting degenerate 2-planes can be $\boldsymbol{c}$-$\boldsymbol{c}$-kind (in 3-dimensional pseudo-Euclidean space isotropic orthogonal vectors do not exist, thus cannot exist $\boldsymbol{c}$-$\boldsymbol{c}$-kind degenerate 2-plane).

**Appendix 3.**

Formula (39) requires the presence, on its right-hand side (RHS), of vectors orthogonal to the common vector $\boldsymbol{d}$ of the intersecting 2-planes $\mathrm{span}[\boldsymbol{u},\boldsymbol{d}]$ and $\mathrm{span}[\boldsymbol{a},\boldsymbol{b}]$ (provided that all these vectors exist and are not isotropic; for isotropic ones see Appendix 2). From such point of view formula (39) (in complex space) has to be presented as

$$\frac{\left|\langle \boldsymbol{u}\wedge\boldsymbol{d}|\boldsymbol{a}\wedge\boldsymbol{b}\rangle\right|^2}{(\boldsymbol{u}\wedge\boldsymbol{d})^2(\boldsymbol{a}\wedge\boldsymbol{b})^2}=\frac{\left|\langle \boldsymbol{u}\wedge\boldsymbol{d}|\boldsymbol{u}_{\|}\wedge\boldsymbol{d}\rangle\right|^2}{(\boldsymbol{u}\wedge\boldsymbol{d})^2(\boldsymbol{u}_{\|}\wedge\boldsymbol{d})^2}=\frac{\left|\langle \boldsymbol{u}|\boldsymbol{u}_{\|}\rangle\right|^2}{\boldsymbol{u}^2\boldsymbol{u}_{\|}^2}, \tag{A3.1}$$

where

$$\boldsymbol{u}_{\|}=\alpha\boldsymbol{a}+\beta\boldsymbol{b}\in\mathrm{span}[\boldsymbol{a},\boldsymbol{b}],\ \langle\boldsymbol{u}_{\|}|\boldsymbol{d}\rangle=0,\ \langle\boldsymbol{u}|\boldsymbol{d}\rangle=0.$$

In order to calculate such $\boldsymbol{u}_{\|}$ let us present

$$\boldsymbol{u}=\alpha\boldsymbol{a}+\beta\boldsymbol{b}+\boldsymbol{u}_{\perp}=\boldsymbol{u}_{\|}+\boldsymbol{u}_{\perp},\ \langle\boldsymbol{a}|\boldsymbol{u}_{\perp}\rangle=\langle\boldsymbol{b}|\boldsymbol{u}_{\perp}\rangle=0. \tag{A3.2}$$

Then, if 2-plane $\mathrm{span}[\boldsymbol{a},\boldsymbol{b}]$ is not degenerate, that is if $D=\boldsymbol{a}^2\boldsymbol{b}^2-\left|\langle\boldsymbol{a}|\boldsymbol{b}\rangle\right|^2=(\boldsymbol{a}\wedge\boldsymbol{b})^2\neq 0,$ one gets

$$\left.\begin{aligned}\langle\boldsymbol{a}|\boldsymbol{u}\rangle&=\alpha\boldsymbol{a}^2+\beta\langle\boldsymbol{a}|\boldsymbol{b}\rangle\\ \langle\boldsymbol{b}|\boldsymbol{u}\rangle&=\alpha\langle\boldsymbol{b}|\boldsymbol{a}\rangle+\beta\boldsymbol{b}^2\end{aligned}\right\}\Rightarrow\quad\begin{aligned}\alpha&=D^{-1}\left(\boldsymbol{b}^2\langle\boldsymbol{a}|\boldsymbol{u}\rangle-\langle\boldsymbol{a}|\boldsymbol{b}\rangle\langle\boldsymbol{b}|\boldsymbol{u}\rangle\right),\\ \beta&=D^{-1}\left(\boldsymbol{a}^2\langle\boldsymbol{b}|\boldsymbol{u}\rangle-\langle\boldsymbol{b}|\boldsymbol{a}\rangle\langle\boldsymbol{a}|\boldsymbol{u}\rangle\right).\end{aligned} \tag{A3.3}$$

Hence,

$$\boldsymbol{u}_{\|}=D^{-1}\left[\boldsymbol{a}\left(\boldsymbol{b}^2\langle\boldsymbol{a}|\boldsymbol{u}\rangle-\langle\boldsymbol{a}|\boldsymbol{b}\rangle\langle\boldsymbol{b}|\boldsymbol{u}\rangle\right)+\boldsymbol{b}\left(\boldsymbol{a}^2\langle\boldsymbol{b}|\boldsymbol{u}\rangle-\langle\boldsymbol{b}|\boldsymbol{a}\rangle\langle\boldsymbol{a}|\boldsymbol{u}\rangle\right)\right]. \tag{A3.4}$$

Evidently, $\langle\boldsymbol{u}|\boldsymbol{u}_{\|}\rangle=\langle\boldsymbol{u}_{\|}+\boldsymbol{u}_{\perp}|\boldsymbol{u}_{\|}\rangle=\boldsymbol{u}_{\|}^2,$ and according to (39), formula (A3.1) gives

$$\frac{\left|\langle \boldsymbol{u}\wedge\boldsymbol{d}|\boldsymbol{a}\wedge\boldsymbol{b}\rangle\right|^2}{(\boldsymbol{u}\wedge\boldsymbol{d})^2(\boldsymbol{a}\wedge\boldsymbol{b})^2}=\frac{\boldsymbol{u}_{\|}^2}{\boldsymbol{u}^2}. \tag{A3.5}$$

Comparing the last formula with expressions in (48) and (47) we find

$$\boldsymbol{u}_{\|}^2=D^{-1}\left(\boldsymbol{a}\langle\boldsymbol{u}|\boldsymbol{b}\rangle-\boldsymbol{b}\langle\boldsymbol{u}|\boldsymbol{a}\rangle\right)^2=D^{-1}|k|^{-2}\boldsymbol{d}^2. \tag{A3.6}$$

Evidently, the same result one gets according to (A3.4).

The following interesting relationships should be noted. According to (18) and (18') we have

$$\frac{(\boldsymbol{a}\wedge\boldsymbol{b}\wedge\boldsymbol{u})^2}{(\boldsymbol{a}\wedge\boldsymbol{b})^2}=\det{}^{-1}\begin{bmatrix}\boldsymbol{a}^2 & \langle\boldsymbol{a}|\boldsymbol{b}\rangle\\ \langle\boldsymbol{b}|\boldsymbol{a}\rangle & \boldsymbol{b}^2\end{bmatrix}\det\begin{bmatrix}\boldsymbol{a}^2 & \langle\boldsymbol{a}|\boldsymbol{b}\rangle & \langle\boldsymbol{a}|\boldsymbol{u}\rangle\\ \langle\boldsymbol{b}|\boldsymbol{a}\rangle & \boldsymbol{b}^2 & \langle\boldsymbol{b}|\boldsymbol{u}\rangle\\ \langle\boldsymbol{u}|\boldsymbol{a}\rangle & \langle\boldsymbol{u}|\boldsymbol{b}\rangle & \boldsymbol{u}^2\end{bmatrix}=$$

$$=\boldsymbol{u}^2+\frac{\langle\boldsymbol{a}|\boldsymbol{u}\rangle}{D}\det\begin{bmatrix}\langle\boldsymbol{b}|\boldsymbol{a}\rangle & \boldsymbol{b}^2\\ \langle\boldsymbol{u}|\boldsymbol{a}\rangle & \langle\boldsymbol{u}|\boldsymbol{b}\rangle\end{bmatrix}-\frac{\langle\boldsymbol{b}|\boldsymbol{u}\rangle}{D}\det\begin{bmatrix}\boldsymbol{a}^2 & \langle\boldsymbol{a}|\boldsymbol{b}\rangle\\ \langle\boldsymbol{u}|\boldsymbol{a}\rangle & \langle\boldsymbol{u}|\boldsymbol{b}\rangle\end{bmatrix}=$$

$$=\boldsymbol{u}^2+\frac{2\mathrm{Re}\langle\boldsymbol{a}|\boldsymbol{u}\rangle\langle\boldsymbol{u}|\boldsymbol{b}\rangle\langle\boldsymbol{b}|\boldsymbol{a}\rangle-\boldsymbol{b}^2\left|\langle\boldsymbol{a}|\boldsymbol{u}\rangle\right|^2-\boldsymbol{a}^2\left|\langle\boldsymbol{b}|\boldsymbol{u}\rangle\right|^2}{D}=\boldsymbol{u}^2-\frac{\left(\boldsymbol{a}\langle\boldsymbol{u}|\boldsymbol{b}\rangle-\boldsymbol{b}\langle\boldsymbol{u}|\boldsymbol{a}\rangle\right)^2}{D}.$$

Therefore, (A3.6) and (A3.5) give

$$\frac{(\boldsymbol{a}\wedge\boldsymbol{b}\wedge\boldsymbol{u})^2}{(\boldsymbol{a}\wedge\boldsymbol{b})^2}=\boldsymbol{u}^2-\boldsymbol{u}_{\|}^2=\boldsymbol{u}_{\perp}^2, \tag{A3.7}$$

$$\frac{\left|\langle \boldsymbol{u}\wedge\boldsymbol{d}|\boldsymbol{a}\wedge\boldsymbol{b}\rangle\right|^2}{(\boldsymbol{u}\wedge\boldsymbol{d})^2(\boldsymbol{a}\wedge\boldsymbol{b})^2}=1-\frac{\boldsymbol{u}_{\perp}^2}{\boldsymbol{u}^2}=1-\frac{(\boldsymbol{a}\wedge\boldsymbol{b}\wedge\boldsymbol{u})^2}{\boldsymbol{u}^2(\boldsymbol{a}\wedge\boldsymbol{b})^2}. \tag{A3.8}$$